\documentclass[12pt]{amsart}
\usepackage{epic,eepic,latexsym, amssymb, amscd, amsfonts, xypic}

\input xy
\xyoption{all}

 \newlength{\baseunit}               
 \newcount{\numlines}                
\newcommand{\point}{\vspace{3mm}\par\refstepcounter{subsection}\noindent{\bf \thesubsection.} }
\newcommand{\tpoint}[1]{\vspace{3mm}\par\refstepcounter{subsection}\noindent{\bf \thesubsection.} 
  {\em #1. ---} }
\newcommand{\epoint}[1]{\vspace{3mm}\par\refstepcounter{subsection}\noindent{\bf \thesubsection.} 
  {\em #1.} }
\newcommand{\bpoint}[1]{\vspace{3mm}\par\refstepcounter{subsection}\noindent{\bf \thesubsection.} 
  {\bf #1.} }

\newcommand{\bpf}{\noindent {\em Proof.  }}
\newcommand{\epf}{\qed \vspace{+10pt}}

\newcommand{\Aut}{\operatorname{Aut}}
\newcommand{\Ext}{\operatorname{Ext}}
\newcommand{\EExt}{\operatorname{\mathcal{E}xt}}

\newcommand{\reldef}{\operatorname{RelDef}}
\newcommand{\relobst}{\operatorname{RelOb}}
\newcommand{\Sing}{\operatorname{Sing}}
\newcommand{\thmstar}{Theorem~$\star$}
\newcommand{\vir}{\text{vir}}
\newcommand{\glued}{\text{glued}}
\newcommand{\conn}{\text{conn}}
\newcommand{\unrigid}{\sim}
\newcommand{\simple}{\text{simple}}
\newcommand{\thick}{\text{thick}}
\newcommand{\thin}{\text{thin}}

\newcommand{\HH}{\mathbb{H}}

\newcommand{\GG}{{\Gamma', \Gamma''}}
\newcommand{\cFbarGG}{{\overline{\mathcal{F}}_{\GG}}}
\newcommand{\cLL}{\mathcal{L}}
\newcommand{\cS}{\mathcal{S}}

\newcommand{\A}{\mathbb{A}}
\newcommand{\Q}{\mathbb{Q}}

\newcommand{\C}{\mathbb{C}}

\newcommand{\proj}{\mathbb P}
\newcommand{\oh}{{\mathcal{O}}}
\newcommand{\cm}{{\mathcal{M}}}
\newcommand{\cmbar}{\overline{\cm}}
\newcommand{\M}{\cmbar}
\newcommand{\mgnbar}{\cmbar_{g,n}}
\newcommand{\mgn}{\mgnbar}
\newcommand{\fM}{{\mathfrak M}}

\newcommand{\tensor}{\otimes}
\newcommand{\cmt}{\widetilde{\cm}}
\newcommand{\mgnrt}{\cm_{g,n}^{rt}} 
\newcommand{\mgnct}{\cm_{g,n}^{ct}} 
\newcommand{\Cech}{\v{C}ech }
\newcommand{\absM}{\left| M \right|}

\newcommand{\cT}{{\mathcal{T}}}
\newcommand{\al}{\alpha}

\newcommand{\be}{\beta}

\newcommand{\de}{\delta}
\newcommand{\De}{\Delta}

\newcommand{\la}{\lambda}

\newcommand{\Hom}{\operatorname{Hom}}

\newcommand{\Spec}{\operatorname{Spec}}

\newcommand{\Spf}{\operatorname{Spf}}

\newcommand{\Bl}{\operatorname{Bl}}
\newcommand{\Sym}{\operatorname{Sym}}
\newcommand{\vdim}{\operatorname{vdim}}

\newcommand{\lremind}[1]{{}}
\newcommand{\bremind}[1]{{}}

\newcommand{\cut}[1]{}

\begin{document}
\pagestyle{plain} \title{{ \large{ Relative virtual localization and
      vanishing of tautological classes on moduli spaces of curves
}
}
}
\author{Tom Graber}
\address{Dept.\ of Mathematics, University of California, Berkeley CA} 
\email{graber@math.berkeley.edu}
\author{Ravi Vakil}
\address{Dept.\ of Mathematics, Stanford University, Stanford CA}
\email{vakil@math.stanford.edu}
\thanks{The first author is partially supported by NSF Grant DMS--0301179
and an Alfred P. Sloan Research Fellowship.  The second author is partially supported by NSF Grant DMS--0238532,  an AMS Centennial Fellowship, and an Alfred P. Sloan Research Fellowship. \newline
}
\date{December 20, 2004.}
\begin{abstract}
  We prove a localization formula for the moduli space of stable
  relative maps.  As an  application, we prove that all
  codimension $i$ tautological classes on the moduli space of stable
  pointed curves vanish away from strata corresponding to curves with at least
  $i-g+1$ genus $0$ components.  As consequences, we prove and
  generalize various conjectures and theorems about various moduli
  spaces of curves (due to Diaz, Faber, Getzler, Ionel, Looijenga,
  Pandharipande, and others). This theorem appears to be the
  geometric content behind these results; the rest is straightforward
  graph combinatorics.  The theorem also suggests the importance of
  the stratification of the moduli space by number of rational
  components.
\end{abstract}
\maketitle
\tableofcontents

{\parskip=12pt 
\section{Introduction}
``Relative virtual localization'' (Theorem~\ref{rvlt}) is a
localization formula for the moduli space of stable relative maps,
using the formalism of \cite{vl}.    It is straightforward to postulate
the form of such a formula. Proving it requires three key ingredients.
\begin{enumerate}
\item[(A)]  The moduli space is shown to admit a $\C^*$-equivariant
locally closed immersion into a smooth Deligne-Mumford stack.  
\item[(B)]  The natural virtual fundamental class
on the $\C^*$-fixed locus is shown to arise from the $\C^*$-fixed
part of the obstruction theory of the moduli space.
\item[(C)] The virtual normal bundle to the fixed locus is
  determined.
\end{enumerate}
The technical hypothesis (A) is presumably unnecessary, but the reader
is warned that it has as yet not been excised from the proof of
virtual localization \cite{vl}.  In verifying it,
we exhibit the moduli space as a global quotient (Sect.\ \ref{quotient}),
which may be independently useful.

We use J.~Li's algebro-geometric definition of this moduli space, and
his description of its obstruction theory \cite{li1, li2}.  We
note the earlier definitions of relative stable maps in the
differentiable category due to A-M.\ Li and Y.  Ruan, and Ionel and
Parker, \cite{liruan, ip1, ip2}, as well as Gathmann's work in the
algebraic category in genus $0$ \cite{andreas}.  We give two proofs of
relative virtual localization, which we hope will give the reader
insight into the technicalities of the moduli space of relative stable
maps and its properties.

The tautological ring of the moduli space of stable $n$-pointed curves
genus $g$ curves $\cmbar_{g,n}$ (and other partial compactifications of
the moduli space of curves $\cm_{g,n}$) is a subring of the Chow ring
(with rational coefficients), containing (informally) ``all of the
classes naturally arising in geometry''. 
As an application of relative virtual localization, we prove the
following result, announced in \cite{notices}.

\tpoint{\thmstar}  {\em Any tautological class of codimension $i$
on $\cmbar_{g,n}$  vanishes on the open set consisting of strata satisfying
$$
\text{\# genus $0$ components} < i-g+1.$$} \label{thmstar} \lremind{thmstar}

(We recall the definition of tautological classes in
Section~\ref{starsection}.)  Equivalently, any tautological class of
codimension $i$ is the pushforward of a class supported on the locus
where the number of genus $0$ components is at least $i-g+1$.  The
corresponding result is not true for Chow classes in general;
$(g,n,i)=(1,11,11)$ provides a counterexample (see \cite[Rem.\
1.1]{socle}).

We emphasize that the proof of \thmstar {} is quite short and naive.
(This portion of the paper 
may be read independently of the others, assuming only the
form of relative virtual localization.)  We define {\em Hurwitz
classes} in the Chow group of the moduli space of curves as the
pushforwards of the (virtual) class of maps to $\proj^1$ with some
branch points fixed.  There are a very small number of tools in the
Gromov-Witten theorist's toolkit, and we employ two of them.  By
deforming the target $\proj^1$, we show that Hurwitz classes lie in a
deep stratum where the curve has many genus $0$ components.  Using
(relative virtual) localization, we express tautological classes as
linear combinations of Hurwitz classes, using the ``polynomiality
trick'' of \cite{socle}.

In Section~\ref{consequences} we give numerous consequences of
\thmstar{}, proving many ``vanishing'' conjectures and theorems on
moduli spaces of curves (due to Diaz, Faber, Getzler, Ionel, Looijenga,
Pandharipande, and others).  This section may also be read
independently of the others.  \thmstar{} implies these various results
via straightforward graph combinatorics, and in most cases extends
them.  The morals of \thmstar{} seem to be: {\em i)} this result is
the fundamental geometry behind these results, and
{\em ii)} the strange stratification of the moduli space of curves by
number of rational components may be potentially important.  In
particular, \thmstar{} raises some natural questions.  For example,
Looijenga asks if this ``$\geq k$-rational component locus'' has the
pseudo-convexity property of being the union of $k+g-1$ affine open sets.
This would imply the Diaz-type results of Section \ref{diazgen} for
example.  See also \cite[Prob.\ 6.5]{hl} for topological and
cohomological vanishing consequences.

We expect relative virtual localization to have many applications;
some of these consequences have already appeared.  For example, as
described in \cite[Sect.~5]{gvelsv}, it gives an immediate proof of
the ``ELSV formula'' \cite{elsv1, elsv2, gvelsv}.  A. Gathmann
has computed the genus $1$ Gromov-Witten invariants 
of the quintic threefold using relative virtual localization \cite{ga2}.
Further applications
are in \cite{fprecent} (including different proofs of some of the
results of Section~\ref{consequences}), and \cite{melissa1,
  melissa2}.  Relative virtual localization also plays an essential
role in Okounkov and Pandharipande's proof of the Virasoro constraints
for target curve \cite{op}.

Throughout this paper we work over $\C$, or any other algebraically
closed field of characteristic $0$.  Our methods are purely algebraic.

\epoint{Acknowledgments} 
We are grateful to J. Li, R. Pandharipande, E.-N.\ Ionel, 
A.~Gathmann, and Y. Ruan for useful
conversations, and to C. Faber for inspiration.  Most of this
work was carried out when the authors were at Harvard University and M.I.T.~
respectively, and it is a pleasure to acknowledge these institutions here.

\section{Stable relative maps}

\label{rvl}\lremind{rvl}In 
order to make this paper as self-contained as possible, we review
relevant background on the moduli space of stable relative maps.  We
give J.~Li's algebro-geometric definition of stable relative maps and
their moduli stack. (We again emphasize the prior work of Li-Ruan and
Ionel-Parker in the differentiable category.)  We then define stable
relative maps to a non-rigid target (implicit in \cite{li1} and
the earlier work in the differentiable category).  We describe the moduli
space of targets $\cT$ and $\cT_{\unrigid}$.  We then show that the
moduli space of stable relative maps satisfies the technical
hypothesis (A) of virtual localization, by exhibiting it as a certain
global quotient.  Finally, we give a description of the perfect
obstruction theory of the moduli space in terms of a relative
obstruction theory over the moduli stack of source and target, which
is a useful description for applications, including the proof of
\thmstar{} in Section~\ref{starsection}.

\bpoint{Notation}\label{notation}\lremind{notation}We first set
notation for the next two sections.  Suppose that $X$ is a smooth
complex projective variety, and $D$ is a smooth divisor.  Let $Y$
be the $\proj^1$-bundle over $D$ corresponding to $N_{D/X} \oplus
\oh_D$: $Y = \proj_D ( N_{D/X} \oplus \oh_D)$.  The bundle $Y$ has two
natural disjoint sections, one with normal bundle $N_{D/X}^\vee$ and the other
with normal bundle $N_{D/X}$; call these the zero section and the
infinity section of $Y$ respectively.  For any non-negative integer
$l$, define $Y_l$ by gluing together $l$ copies of $Y$, where the
infinity section of the $i^{\text{th}}$ component is glued to the zero section
of the $(i+1)^{\text{st}}$ ($1 \leq i < l$).  Denote the zero section of the
$i^{\text{th}}$ component by $D_{i-1}$, and the infinity section by $D_i$, so
$\Sing Y_l = \cup_{i=1}^{l-1} D_i$.  We will also denote $D_l$ by
$D_{\infty}$.  Define $X_l$ by gluing $X$ along $D$ to $Y_l$ along
$D_0$.  (Thus for example $\Sing X_l = \cup_{i=0}^{l-1} D_i$, and $X_0=X$.)
Let $\Aut_D Y_l\cong (\C^*)^l$ be the obvious group of
automorphisms of $Y_l$ preserving $D_0$, $D_\infty$, and 
the morphism to $D$
 and let $\Aut_D
X_l$ be the group of automorphisms of $X_l$ preserving $X$ (and $D$) and with
restriction to $Y_l$ contained in $\Aut_D Y_l$ (so $\Aut_D X_l \equiv \Aut_D Y_l
\cong (\C^*)^l$).  Note that if $D$ is disconnected, these conventions differ
from the most natural interpretation of the notation which would be strictly
larger groups.  This will be irrelevant for the applications in the sequel, but
important in general.

\bpoint{The stack of stable relative maps} We recall the definition of
stable relative maps to $X$ relative to a divisor $D$.  If $\Gamma$ denotes the data of
\begin{enumerate}
\item[A1.] arithmetic genus $g$ of the source curve, 
\item[A2.] element $\be$ of $H_2(X)$,
\item[A3.] number $n$ of marked points mapped to $D$, and corresponding
partition $\alpha$ of $\be \cdot D$ into $n$ parts, 
$\alpha_1, \ldots, \alpha_n$, and 
\item[A4.] number $m$ of other marked points,
\end{enumerate}
we denote the moduli stack by
$\cmbar_{\Gamma}(X,D)$. 

We make two conventions about stable relative maps to $X$ relative to
a divisor $D$ differing from Li's: the source curve need not be
connected, and the points mapping to $D$ are labeled.  These
conventions are not mathematically important, but simplify the
exposition.  (Hence for us, moduli spaces of stable relative maps are
not in general connected, and the arithmetic genus of the source curve
may be negative.)

The $\C$-points of this stack correspond to morphisms $C \stackrel f
\rightarrow X_l  \rightarrow X$ where $C$ is a nodal
curve of arithmetic genus $g$ equipped with a set of marked smooth points,
$p_1, \ldots , p_m$, $q_1, \ldots, q_n$.  The morphism $f$ has the
property that $f^*(D_{\infty})=\sum \alpha_i q_i$, and it satisfies
the {\em predeformability condition} above the singular locus
$\Sing X_l = \cup_{i=0}^{l-1} D_i$ of $X_l$, meaning that the
preimage of the singular locus is a union of nodes of $C$,
and if $p$ is one such node, then the two branches of $C$ at $p$ map
into different irreducible components of $X_l$, and their orders of
contact with the divisor $D_i$ (in their respective components of
$X_l$) are equal.  We will refer to these nodes as
the {\em distinguished nodes} of $C$.  
The morphism $f$ is also required to satisfy a
stability condition that there are no infinitesimal automorphisms of
the sequence of maps $(C, p_1, \dots, p_m, q_1, \dots, q_n)
\rightarrow X_l \rightarrow X$, where the allowed automorphisms of the
map from $X_l$ to $X$ are $\Aut_D(X_l)$.

\point \lremind{plug}\label{plug}The 
paper \cite{li1} defines a good notion of a family of such maps,
i.e.\ a moduli functor or groupoid.  A family of stable relative maps
over a base scheme $S$ is a pair of morphisms of flat $S$-schemes
$C \stackrel f \rightarrow \overline X \rightarrow X \times S$, where for each
$\C$-point $s$ in $S$, the fiber $C_s \stackrel {f_s} \rightarrow \overline X_s
\rightarrow X$ gives a stable relative map.  There is also the
predeformability condition, that in a neighborhood of a node of $C_s$
mapping to a singularity of $\overline X_s$, we can choose
\'etale-local coordinates on $S$, $C$, and $\overline X$ with 
charts  of the form $\Spec R$, 
$\Spec R[u,v]/uv=b$ 
and 
$\Spec R[x,y,z_1,\ldots , z_k]/xy=a$ 
respectively, 
with the map of
the form $x \mapsto \alpha u^m$, $y \mapsto \beta v^m $ with
$\alpha$ and $\beta$ units, and no
restriction on the $z_i$.  Given any family of stable morphisms to an
allowed family of target schemes,
the locus of maps satisfying this predeformability condition naturally
forms a locally closed subscheme.

\bpoint{Stable relative maps to a  non-rigid target} 
\label{nonrigid}\lremind{nonrigid}We 
make explicit the following variation on stable relative maps.
Such maps appear in the boundary of the space of ``usual'' stable
relative maps \cite[Sect.~3.1]{li1}.  (For this definition, $Y$ can be
any $\proj^1$-bundle with two disjoint sections 
over a smooth variety $D$, but we will be
using this definition in the context described in
Section~\ref{notation}.)

We consider morphisms to $(Y_l, D_0 \cup D_{\infty})$, where
two stable maps 
$$\xymatrix {(C, p_1, \dots, p_m, q_1, \dots, q_n, r_1, \dots,
r_{n'}) \ar[r] &  (Y_l, D_0, D_{\infty})}$$
($r_i \rightarrow D_0$, $q_i
\rightarrow D_{\infty}$) and 
$$\xymatrix{  (C', p'_1, \dots, p'_m, q'_1, \dots, q'_n,
r'_1, \dots, r'_{n'}) \ar[r] & (Y_l, D_0, D_{\infty})}$$
are considered
isomorphic if there is a commutative diagram
$$
\xymatrix{ (C, p_1, \dots, r_{n'}) \ar[r]^{\sim} \ar[d] & 
    (C', p'_1, \dots, r'_{n'}) \ar[d] \\
(Y_l, D_0, D_{\infty}) \ar[r]^{\alpha}  & (Y_l, D_0, D_{\infty})}
$$
with $\alpha \in \Aut_D Y_l$.  Predeformability
is required above the singularities of $Y_l$.    
We will call these {\em stable relative
maps to a non-rigid target}.
The definition of a family of such maps is the obvious 
variation of that for stable relative maps described in the
previous section
(see also \cite[Sect.~2.2]{li1}).

Such maps also form a Deligne-Mumford moduli
stack. 
If $\Gamma$ denotes
the data of 
\begin{enumerate}
\item[B1.] arithmetic genus of the source curve, 
\item[B2.] element $\be$ of $H_2(Y)$,
\item[B3.] number of marked points mapped to $D_0$ (resp.\ $D_{\infty}$), and corresponding
partition of $\be \cdot D_0$  (resp.\ $D_{\infty}$), and
\item[B4.] number of other marked points,  
\end{enumerate}
we denote the corresponding moduli stack
$\cmbar_{\Gamma}(Y, D_0, D_{\infty})_{\unrigid}$.  This stack
may be constructed by a straightforward variation of the
global quotient construction of Section \ref{quotient}.


\bpoint{The moduli spaces of targets $\cT$ and $\cT_{\unrigid}$}
 \lremind{modtarget}\label{modtarget}\lremind{psi}\label{psi}Let $\cT$
 be the Artin stack parametrizing the possible targets of
 relative stable maps to $(X,D)$.  (This is called the stack of
 expanded degenerations in \cite[Sect.\ 1]{li1}.)  It has one
 $\C$-point for each nonnegative integer.  $\cT$ is isomorphic to the
 open substack of the Artin stack $\fM_{0,3}$ of prestable
 $3$-pointed genus $0$ curves consisting of curves where the only
 nodes separate point $\infty$ from points $0$ and $1$ (see Figure~\ref{tt}).  
 On the level of $\C$-points, this open immersion corresponds to
 replacing $X_l$ with a copy of $\proj^1$ (with $0$ and $1$ marked)
 attached at $\infty$ to the fiber of $Y_l$ over a fixed point of $D$,
 marking the point on $D_\infty$.  The equality of these stacks can be
 seen either by describing this construction for families (using a
 blow-up construction rather than gluing), or by noting
 that the construction in \cite{li1} is independent of $(X,D)$.  Note
 in particular that $\cT$ is independent of $(X,D)$, but the universal
 family is not.  We refer the reader to \cite{li1} for the
 construction of the universal family over $\cT$.

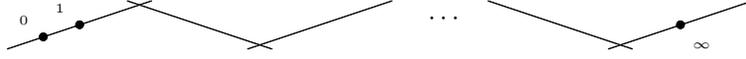
\begin{figure}
\begin{center}
\setlength{\unitlength}{0.00083333in}
\begingroup\makeatletter\ifx\SetFigFont\undefined%
\gdef\SetFigFont#1#2#3#4#5{%
  \reset@font\fontsize{#1}{#2pt}%
  \fontfamily{#3}\fontseries{#4}\fontshape{#5}%
  \selectfont}%
\fi\endgroup%
{\renewcommand{\dashlinestretch}{30}
\begin{picture}(4674,363)(0,-10)
\put(4212,186){\blacken\ellipse{50}{50}}
\put(4212,186){\ellipse{50}{50}}
\path(3012,336)(3912,36)
\path(3012,336)(3912,36)
\path(3762,36)(4662,336)
\path(3762,36)(4662,336)
\put(4287,36){\makebox(0,0)[lb]{\smash{{{\SetFigFont{5}{6.0}{\rmdefault}{\mddefault}{\updefault}$\infty$}}}}}
\put(237,111){\blacken\ellipse{50}{50}}
\put(237,111){\ellipse{50}{50}}
\put(462,186){\blacken\ellipse{50}{50}}
\put(462,186){\ellipse{50}{50}}
\path(12,36)(912,336)
\path(12,36)(912,336)
\path(762,336)(1662,36)
\path(762,336)(1662,36)
\path(1512,36)(2412,336)
\path(1512,36)(2412,336)
\put(312,261){\makebox(0,0)[lb]{\smash{{{\SetFigFont{5}{6.0}{\rmdefault}{\mddefault}{\updefault}$1$}}}}}
\put(2637,186){\makebox(0,0)[lb]{\smash{{{\SetFigFont{10}{12.0}{\rmdefault}{\mddefault}{\updefault}$\cdots$}}}}}
\put(87,186){\makebox(0,0)[lb]{\smash{{{\SetFigFont{5}{6.0}{\rmdefault}{\mddefault}{\updefault}$0$}}}}}
\end{picture}
}
\end{center}
\caption{The open substack $\cT$ of $\overline{\mathfrak{M}}_{0,3}$:
all nodes separate $\infty$ from $0$ and $1$\lremind{tt}}
\label{tt}
\end{figure}

The moduli space $\M_\Gamma(X,D)$ is then simply the locally closed
 subset of the space of stable morphisms to the universal family over
 $\cT$ where the predeformability and stability conditions are
 satisfied.

The analogous space for stable relative maps to non-rigid targets
is the stack $\cT_\unrigid$ which is the open substack of $\fM_{0,2}$
parametrizing $2$-pointed curves where the only nodes separate $0$
from $\infty$.  The line bundle corresponding to the cotangent space
at the point $0$ will play an important role; we denote its first
Chern class by $\psi$.  (The notation $\psi_0$ is more usual, but we
wish to avoid confusing this class with the $\psi$-classes coming from
the marked points on the source of the relative stable map.)  We
denote the pullback of $\psi$ to $\M_{\Gamma}(Y, D_0, D_{\infty})_\unrigid$ by
$\psi$ as well.  See \cite{eric} for a different definition of $\psi$,
and a more detailed discussion.

\bpoint{Global equivariant embedding, and $\cmbar_{\Gamma}(X,D)$ 
as a global quotient} \label{quotient}\lremind{quotient}To apply the virtual
localization theorem of \cite{vl} we need to verify
the technical hypothesis that the moduli space $\cmbar_{\Gamma}(X,D)$
admits a $\C ^*$-equivariant locally closed immersion in a smooth
Deligne-Mumford stack (in the case where $(X,D)$ admits a $\C^*$-action). 
In \cite{vl}, this is verified for the
space of ordinary stable maps 
to a smooth projective variety with $\C^*$-action.
The method used there is to realize the moduli space as a quotient stack
$[V/G]$ with $V$ a quasi-projective variety and $G$ a reductive group
such that the following
two conditions are satisfied:
\begin{itemize}
\item $V$ is a locally closed subset of a smooth projective variety
$W$ such that the $G$-action on $V$ extends to an action on $W$.
\item There is a $\left( \C^* \times G \right)$-action on $W$ which preserves $V$
and descends to the $\C^*$-action  on the moduli stack.
\end{itemize}
We  give an analogous 
construction for the space of relative stable 
maps.  
To do this, we use the constructions of \cite{fup} and
\cite{vl} for stable maps.
In \cite[Sect.\ 2]{fup}, the moduli space of stable maps
$\cmbar_{g,n}(\proj, \be)$ to a projective space $\proj$ (with fixed
numerical data $g$, $n$, $\be$) is expressed as a quotient of a
locally closed subscheme $J$ of a Hilbert scheme $H$ (of a product of
projective spaces) by a reductive group $G$.
As observed in \cite[App.\ A]{vl}:
\begin{enumerate}
\item[(i)] The stack quotient $[J / G]$ is $\cmbar_{g,n}(\proj, \be).$
\item[(ii)] There is a $(\C^* \times G)$-action on $J$ which descends
to the given $\C^*$-action on $\cmbar_{g,n}(\proj, \be)$.
\item[(iii)]  There is a $(\C^* \times G)$-equivariant 
linearized locally closed immersion of $J \subset H$ in a (smooth) Grassmannian
$\mathbb{G}$.
\end{enumerate}

Now fix data $\Gamma$ satisfying A1--A4.  For any stable relative 
map in
$\cmbar_{\Gamma}(X,D)$ 
$$\xymatrix{f: C \ar[r] &  X_l}$$ 
where the target breaks
(i.e.\ $l \geq 1$), and any choice of  singular
locus $D_i$ of $X_l$, let $\Gamma'$ and $\Gamma''$ be the ``splitting data''
of the data $\Gamma$ that arises by separating $X_l$ into two pieces.
(Thus the stable relative map can be interpreted as a point of
$\cmbar_{\Gamma'}(X,D) \times_{D^n} \cmbar_{\Gamma''}(Y, D_0,
D_{\infty})_\unrigid$.)  

Let $M$ be the set of possible splitting data $(\GG)$.
Then $M$ is finite, and has a partial ordering defined as
follows: if $f: C \rightarrow X_l$ is a stable relative map where $l
\geq 2$, and breaking the target into $X_{l_1} \cup Y_{l-l_1}$ yields
splitting data $(\Gamma'_1, \Gamma''_1)$, and breaking the target into
$X_{l_2} \cup Y_{l-l_2}$ yields splitting data $(\Gamma'_2,
\Gamma''_2)$, and $0 \leq l_1 < l_2 < l$, then $(\Gamma'_1,
\Gamma''_1) \prec (\Gamma'_2, \Gamma''_2)$.  (The reader should verify
that this is indeed a partial order.)  Extend the partial order to a
total order, i.e.\ choose an order-preserving bijection $\mathbf{i}:
\{ 1, \dots, \absM \} \rightarrow M$.

In $X \times \left(\proj^1 \right)^{\absM}$, define the locus $L_i$
($1 \leq i \leq \absM$) by $D \times \proj^1 \times \cdots \times \{ 0
\} \times \cdots \times \proj^1,$ where the $\{ 0 \}$ appears in the
$i^{\text{th}}$ $\proj^1$-factor.  Blow up $X \times \left(\proj^1
\right)^{\absM}$ along $L_1$, then along the proper transform of
$L_2$, and so on, ending with the proper transform of $L_{\absM}$.
Let $\Bl$ denote the result.  The singular locus of the morphism $\Phi: \Bl
\rightarrow \left(\proj^1 \right)^{\absM}$ is a union of
$\absM$ codimension $2$ loci, say $\Sing(1)$, \dots, $\Sing({\absM})$,
in order-preserving bijection with $\{1, \dots, \absM \}$.  (In other
words, given a fiber of the morphism $\Phi$ isomorphic to $X_l$, 
such that $D_{i_1} \subset \Sing(j_1)$  and
$D_{i_2} \subset \Sing(j_2)$, then $i_1 < i_2 \Leftrightarrow j_1 < j_2$.)

There is a
$(\C^*)^{ 1+\absM}$-action on
  $$\left( X \times \left( \proj^1 \right)^{\absM}, L_1, \dots,
    L_{\absM} \right)$$
(where the $i^{\text{th}}$ $\C^*$ acts on the $i^{\text{th}}$ factor
of $X \times \left( \proj^1 \right)^{\absM}$)
  and hence on $\Bl$.  Choose a
  $(\C^*)^{1+\absM}$-equivariant closed immersion of $\Bl$ into a large
  projective space $\proj$.
(As $X$ is a smooth projective variety with
$\left( \C^* \right)^{1+\absM}$-action, there is a 
$\left( \C^* \right)^{1+\absM}$-equivariant ample invertible sheaf
$\cLL$ on $X$ \cite[Cor.\ 1.6, p.~35]{git}.)

  Then consider the moduli space of stable maps $\cmbar_{g,m+n}(\proj,
  \be')$, where $\be'$ is the homology class of $\proj$ corresponding
  to the curve mapping to $\be \in X$.  
As stated above, this space can be expressed as a
  quotient satisfying (i)--(iii).  Consider the locally closed substack
  $\cmbar'$ of $\cmbar_{g,m+n}(\proj,
  \be')$ corresponding to stable maps:
\begin{enumerate}
\item[(I)] mapping to $\Phi^{-1} 
\left( \left( \A^1\right)^{\absM}\right) \subset \Bl$;
\item[(II)] predeformable as maps to $\Bl \rightarrow 
\left( \proj^1 \right)^{\absM}$;
\item[(III)] for any stable map with image meeting $\Sing(l)$, the
splitting type above that singularity is of type $\mathbf{i}(l)$; and
\item[(IV)] with the pullback of $D_{\infty}$ given by $\sum \al_i q_i$.
\end{enumerate}

Then it is straightforward to verify that the quotient $[\cmbar' / (
\C^*)^{\absM}]$ is isomorphic to $\cmbar_{\Gamma}(X,D)$.  (For
example, from the construction it is clear how to make smooth-local
sections of this quotient map.)  Note that condition (III) is
essential.

Finally, the locally closed immersion $\cmbar' \hookrightarrow
\cmbar_{g,m+n}(\proj, \be')$ lifts to a
locally closed immersion $J' \hookrightarrow J$, and the 
$(\C^*)^{\absM}$-action
lifts to the Grassmannian $\mathbb{G}$.
Hence we have shown the following technical condition required to
apply virtual localization.

\tpoint{Theorem}  {\em The stack $\cmbar_{\Gamma}(X,D)$ admits
a $\C^*$-equivariant closed immersion into a smooth Deligne-Mumford
stack.}

Specifically, we have shown that conditions (i)--(iii) above hold, with
$\cmbar_{g,n}(X,\be)$ replaced by $\cmbar_{\Gamma}(X,D)$; $J$ replaced
by $J'$; and $G$ replaced by $G \times
(\C^*)^{\absM}$.
As a bonus, we observe that 
$\cmbar_{\Gamma}(X,D)$ is a global quotient.

\bpoint{The perfect obstruction theory} A
perfect obstruction theory on $\M_\Gamma(X,D)$ is
constructed and studied
in \cite{li2}.  There it is made explicit that this
obstruction theory is induced by a relative perfect obstruction theory,
relative to the morphism $\M_\Gamma(X,D) \rightarrow \cT$ to the moduli
space of targets (Sect.\ \ref{modtarget}).  In fact, 
an analysis of the obstruction theory constructed there shows that it 
is induced by
a relative obstruction theory over the product of the (smooth) stacks
parametrizing the possible targets and possible sources.  
That is, if
we let $E \rightarrow L_{\M_\Gamma(X,D)}$ be the perfect obstruction
theory on $\M_\Gamma(X,D)$ and we consider the natural morphism
$$ \xymatrix{\Phi : \M_\Gamma(X,D) \ar[r] & \cT \times \fM_{g,m+n}
},$$ then there exists an element $F$ in the derived category of
sheaves on $\M_\Gamma(X,D)$ locally representable by a two-term
complex of vector bundles, a morphism from $F$ to the relative
cotangent complex of the morphism $\Phi$ satisfying the usual
cohomological conditions (see \cite[Sect.\ 4]{bf}), and a distinguished
triangle
$$ \xymatrix{ F[-1] \ar[r] &  \Phi^* L_{\cT \times \fM_{g,m+n}} \ar[r] &
 E  \ar[r] & F } $$
compatible with the morphisms from $E$ and $F$ to the appropriate
cotangent complexes.

The reason this is useful is that $F$ can be understood 
easily in terms of explicit cohomology groups.  One should think
of the cohomology sheaves of $F$ (or its dual) 
as measuring the deformations and
obstructions of a stable relative map, 
once deformations of the source and target are chosen.  We will explain
this carefully here, since this material is spread
through a large portion of \cite{li2} and is not
in the form that we need.

We first introduce some notation.  If $V$ is
a variety, and $D$ is 
a normal crossings divisor contained in the smooth locus of $V$,
then let $\Omega_V (\log D)$ be the sheaf of 
K\"ahler differentials with logarithmic poles along $D$.  
Denote the dual of this sheaf
by $T_V (-\log D)$.
When $V$ has normal crossings singularities, this
can be interpreted  as vector fields on the normalization of $V$ which are
tangent to the divisor $D$, tangent to the singular strata, and such that
the induced vector fields on the singular strata agree on the different
 branches
of the normalization.  When
$V$ is a nodal curve, $T_V(- \log D)$ corresponds to
vector fields vanishing at the divisor $D$ and the nodes of $V$.

Finally, if $f: C \to X_l \to X$ is a relative stable map,
then the sheaf $f^*T_{X_l}(-\log D_{\infty})$ 
will have a torsion subsheaf supported
at the distinguished nodes of $C$.  This sheaf is the kernel
of the morphism from $f^*T_{X_l}(-\log D_\infty)$ to 
$\Hom(f^*\Omega_{X_l}(\log D_\infty , \oh_C)$.  We will use
the notation $f^\dagger T_{X_l}(-\log D_\infty)$ to denote the
quotient of $f^*T_{X_l}(-\log D_\infty)$ by this subsheaf.  Sections
of $f^\dagger T_{X_l}(-\log D_\infty)$ when restricted to an irreducible
component $C'$ of $C$ correspond simply to sections of the pullback of
the restriction of the logarithmic tangent bundle to the irreducible 
component of $X_l$ to which $C'$ maps.  At a distinguished node, the sections
on the two components of $C$ are required to agree in the sense that they give
the same element of the tangent space to $D$.

If $C \stackrel f \rightarrow X_l \rightarrow X$ is a stable relative map 
to $(X,D)$, then it is easy to compute that
the $\Phi$-relative tangent space to
the moduli space at this point is\lremind{postit}
\begin{equation}\label{postit}
\reldef(f) = H^0(C,f^\dagger T_{X_l} (-\log D_{\infty})).
\end{equation}
 It is important not to
confuse this 
space with the space $\Hom(f^* \Omega_{X_l} (\log D_{\infty}),\oh_C)$
which parametrizes deformations of the map preserving the contact
conditions along $D_{\infty}$, but ignoring the 
predeformability condition at the 
nodes. 
  The space of relative obstructions
$\relobst (f)$
has a natural filtration 
\begin{equation}\label{localglobal} 
0 \rightarrow H^1(C,f^*T_{X_l}(-\log D_{\infty})) \rightarrow \relobst (f) 
\rightarrow
H^0(C,f^{-1}\EExt^1(\Omega_{X_l}(\log D_{\infty}), \oh_{X_l})) \rightarrow 0. 
\end{equation}
Because $D_{\infty}$ is contained in the smooth locus of $X_l$, 
there is a natural isomorphism between 
$\EExt^1(\Omega_{X_l}(\log D_{\infty}),\oh_{X_l})$ and 
$\EExt^1(\Omega_{X_l},\oh_{X_l})$.  We can make this sheaf
more explicit as follows.  $\EExt^1(\Omega_{X_l},\oh_{X_l})$
is the pushforward of a trivial line bundle on $\Sing X_l  = \cup_{i=0}^{l-1}D_i$.  
However, there is no canonical trivialization of this line bundle,
and in families it can vary.  On
$D_i$ this line
bundle is canonically isomorphic to  $N_{D_i / X({i-1})} \tensor N_{D_i / X(i)}$,
where $X(i-1)$ and $X(i)$ are the components of $X_l$ containing $D_i$
(we will not use this notation $X(i)$ again).

Sequence \eqref{localglobal} 
is a close analog of the local-to-global sequence for $\Ext$
and has the same deformation-theoretic interpretation:
the right-hand group represents the local obstructions to
deforming $f$, and the left-hand group is the global obstruction.
We now explain explicitly why these groups are as above.  
 
For each singular locus $D_i$ of $X_l$, 
we denote the 
set of nodes of $C$ mapping to $D_i$ by
$\{N^1_i,\ldots , N_i^{j_i}\}$.
We can rewrite sequence \eqref{localglobal} as\lremind{explocalglobal}
\begin{equation} \label{explocalglobal} 
\xymatrix{ 0 \ar[r] &  H^1(C,f^\dagger T_{X_l}(-\log D_{\infty})) \ar[r] & 
\relobst (f) 
\ar[r] & \bigoplus_{i=0}^{l-1} L_i^{\oplus j_i} \ar[r]&  0} 
\end{equation}
where $L_i$ is the one-dimensional vector space of sections of the
sheaf $\EExt^1(\Omega_{X_l}, \oh_{X_l})|_{D_i}$, which can be
interpreted as the deformation space of the singularity $D_i$ of
$X_l$.  The last term is the local obstruction to extending a map,
which comes from a compatibility requirement between the choices of
deformation of a singularity $D_i$ of the target and the deformation
of those nodes of the source $N^j_i$ which map to $D_i$.  The fact
that this obstruction space is identified with the deformation space
of the target singularity $D_i$ can be seen from a direct local calculation.

Explicitly, suppose we have a predeformable morphism between
nodal curves over a base $\Spec R$ with $R$ an Artin local ring
with maximal ideal $\mathfrak m$.  For notational convenience (to avoid
the irrelevant variables $z_i$ of Sect.\ \ref{plug}) 
we assume $\dim X=1$.  \'Etale-locally,
the predeformable morphism has the form $C \rightarrow X_l$ 
with $X_l=\Spec (R[x,y]/(xy-a))$ and
$C= \Spec (R[u,v]/(uv-b))$ with $a,b \in \mathfrak m$
and the morphism is given by
$x \mapsto \alpha u^n, y\mapsto \beta v^n$.  If $\tilde R$ is
a small extension of $R$ with ideal $I$, then choosing
lifts of $C$ and $X_l$ to $\Spec \tilde R$ corresponds to lifting
the elements $a$ and $b$ in $R$ to elements $\tilde a$ and
$\tilde b$ in $\tilde R.$  These choices
are torsors for the ideal $I$, but the structure of $I$-torsor
depends on the choice of local coordinates on $C$ and $X_l$.  
The choice of $\tilde a$ is canonically a torsor for
$I \tensor L_i$ where $L_i$ is the one-dimensional deformation
space of the singularity $D_i$ of $X_l$.

If we want the 
morphism to extend, then the choice of $\tilde b$ determines
the choice of $\tilde a$, because we have the formula
$a=xy= \alpha \beta u^n v^n = \alpha \beta b^n.$  If we choose
a lifting $\tilde b$ of $b$ and then try to choose liftings
$\tilde \alpha, \tilde \beta, \tilde a$ satisfying the
analogous formula, we see that $\tilde a$ is determined,
since the ambiguity in the choice of $\tilde \alpha$ and 
$\tilde \beta$ is an element of $I$ which is annihilated by
multiplication with $\tilde b$.  Hence, given choices of
$\tilde b$ and $\tilde a$, the element 
$\tilde a - \alpha \beta \tilde b^n$ is a canonically defined
element of $I\tensor L_i$ which vanishes exactly if there
is a local extension of the predeformable morphism.  

The first term in Sequence (\ref{localglobal}) is the global obstruction,
analogous to the term $H^1(C,f^*T_X)$ in the case of ordinary stable maps,
which can be identified by the usual \Cech cocycle construction
once the local obstructions vanish.  Explicitly, suppose we are given
a morphism $C \rightarrow X_l$ over $\Spec R$ and liftings $\tilde C$
and $\tilde {X_l}$ as above.  Then if all the local obstructions
constructed above vanish, we can \'etale-locally extend the morphism.
That is, we can choose an \'etale cover of $C$ by open sets $U_i$ 
inducing an \'etale cover of $\tilde C$ by open
sets $\tilde U_i$, and on each $\tilde U_i$ there exists an extension of the
morphism.  The set of such extensions is a 
$H^0(U_i, f^\dagger T_{X_l}(-\log D))$-torsor, 
so taking the differences on overlaps
we get a \Cech cocycle representing an element of 
$H^1(C,f^\dagger T_{X_l}(-\log D))$
which gives the global obstruction. (Here we use the fact
that \'etale cohomology agrees with Zariski cohomology for
coherent sheaves.)

\epoint{The perfect obstruction theory of $\cmbar_{\Gamma}(Y, D_0, D_{\infty})_{\unrigid}$}   
\lremind{pot}\label{pot}The relative obstruction theory for 
$\cmbar_{\Gamma}(Y, D_0, D_{\infty})_{\unrigid}$ (over the moduli stack
of source and target) is identical to
the previous (rigid) case, with $\cT$ replaced by $\cT_{\unrigid}$.

\section{Relative virtual localization} 
\label{svl}\lremind{svl}In \cite{vl}, 
a localization formula is proved for the virtual
fundamental class in the general context of $\C^*$-equivariant perfect
obstruction theories.  If $\cm$ is a Deligne-Mumford stack with such
an obstruction theory, then the $\C^*$-fixed loci $\cm_i$ carry an
associated ($\C^*$-fixed) perfect obstruction theory, giving them each
a virtual fundamental class.  The {\em virtual normal bundle} to
$\cm_i$ arises from the ``moving'' (non-torus-fixed) part of the obstruction theory.
Assume that $\cm$ has an equivariant locally closed immersion into a
smooth Deligne-Mumford stack.  The virtual localization formula
\cite[equ.~(1)]{vl} then states:
$$
[\cm]^{\vir} = \sum \frac { [\cm_i]^{\vir}} { e(N_i^{\vir})}$$
in
$A_*^{\C^*}(\cm) \otimes \Q[t, \frac 1 t]$, where $t$ is the generator
of the $\C^*$-equivariant Chow ring of a point.  (In applications, the
equation is usually capped with the Chow classes of equivariant vector
bundles to obtain equalities of intersection numbers.  Our strategy is to
cap with classes of lower codimension, to prove equalities in the Chow ring.)

The primary example there was the moduli space of stable maps to
$\proj^n$ \cite[Sect.~4]{vl}.  In that case the virtual normal bundle
is given explicitly by contributions from parts of the geometry of the
maps in $\cm_i$.

In this section, we suppose $X$ is a smooth projective variety with a
non-trivial $\C^*$-action, and $D$ is an irreducible divisor in the
fixed locus.  (The arguments below apply essentially without change in
greater generality, for example when $D$ is reducible or the group is
larger, but for the sake of exposition we will state relative virtual
localization in only moderate generality.)  The natural $\C^*$-action
on the perfect obstruction theory of $\cmbar_{\Gamma}(X,D)$ gives this
space an equivariant virtual fundamental class.
 
\epoint{$\C^*$-fixed loci in the relative setting, and their virtual
  fundamental classes}
We now set notation for  the $\C^*$-fixed locus of $\cmbar_{\Gamma}(X,D)$.  If
a connected component of the fixed locus has general morphism with
target $X$ (i.e.\ the target doesn't degenerate for the general
morphism, or equivalently for {\em any} morphism in this connected
component of the fixed locus), we say that it is a
{\em simple} fixed locus.  Otherwise, it is a {\em composite} fixed
locus.

A simple fixed locus has an induced
obstruction theory that is the fixed part of the obstruction theory of
$\cmbar_{\Gamma}(X,D)$, and hence a virtual
fundamental class, and a virtual normal bundle, which we denote
$N_{\Gamma}$.  The analysis of
the perfect obstruction theory here is much simpler than the general case,
because the terms corresponding to deformations and automorphisms of the target
as well as the term giving the local obstruction to deforming the morphism vanish.
Thus the analysis of the virtual normal bundle to such a component is
identical to the case of ordinary stable maps with
the bundle $T_X$ systematically replaced
by $T_X(-\log D)$.  Denote the union of simple fixed loci by
$\cmbar_{\Gamma}(X,D)^{\simple}$.  

Any element of a composite fixed locus is of the form $f: C' \cup C''
\rightarrow X_l$ ($l>0$), which restricts to $f': C' \rightarrow X$
and $f'': C'' \rightarrow Y_l$ which agree over the nodes $\{ N_1,
\dots, N_{\de}\} = C' \cap C''$. Let $\Gamma'$ be the data A1--A4
corresponding to $f'$, and $\Gamma''$ be the data B1--B4 corresponding
to $f''$.  (Any two of $\{ \Gamma, \GG \}$ determine the
third.)  Define $m_i$ by $(f')^{-1}(D) = \sum m_i N_i$ on $C'$, or
equivalently $(f'')^{-1}(D_0) = \sum m_i N_i$ on $C''$; this is part
of the data of both $\Gamma'$ and $\Gamma''$.  Both $\Gamma'$ and
$\Gamma''$ are locally constant on the fixed locus.

The fixed locus $\cFbarGG$ corresponding to a given $\Gamma'$ and
$\Gamma''$ is canonically the (\'etale) quotient of the stack
$$
\cmbar_{\GG} = \cmbar_{\Gamma'}(X,D)^{\simple} \times_{D^n}
\cmbar_{\Gamma''}(Y,D_0, D_{\infty})_{\unrigid},
$$
by the finite group $\Aut(m_\cdot)$ (those permutations of $(1,
\dots, \de)$ preserving $(m_1 \dots, m_\de)$).  Call 
this quotient morphism $gl$, so
$$
gl: \cmbar_{\GG} \rightarrow \cFbarGG.$$
(Cf.\ \cite[Prop.\ 4.13]{li1}; the
morphism $gl$ is called $\Phi$ there.)

A virtual fundamental class on $\cmbar_{\GG}$ (and $\cFbarGG$), which we
term the {\em glued} fundamental class, is induced by the virtual
fundamental classes on the factors (cf.\ \cite[p.~203]{li2}):\lremind{butter}
\begin{eqnarray} \label{butter}
[\cmbar_{\Gamma',
  \Gamma''}]^{\glued} &=& \Delta^!
\left(
[\cmbar_{\Gamma'}(X,D)^{\simple}]^{\vir} \times
[\cmbar_{\Gamma''}(Y,D_0, D_{\infty})_{\unrigid}]^{\vir} \right)  \\
\nonumber
{}[ \cFbarGG ]^{\glued}
&=& \frac 1 { \left| \Aut( m_{\cdot}) \right|}
gl_* [\cmbar_{\Gamma',
  \Gamma''}]^{\glued} \end{eqnarray}
where $\Delta: D^n \rightarrow D^n \times D^n$ is the diagonal
morphism.  A second virtual fundamental class on $\cFbarGG$ (and by
pullback, on $\cmbar_{\Gamma',
  \Gamma''}$), $[\cFbarGG]^{\vir}$ (resp.\ $[\cmbar_{\GG}]^{\vir}$), 
is that induced
by the $\C^*$-fixed part of the pullback of the obstruction theory of
$\cmbar_{\Gamma}(X,D)$.  The following lemma shows 
that they are the same.

\tpoint{Lemma} {\em \lremind{butterfly}\label{butterfly}The 
torus-induced virtual fundamental class on a composite
fixed locus is exactly the glued virtual fundamental class coming
from the two factors.  In other words, $$[\cFbarGG]^{\glued} =
[\cFbarGG]^{\vir} \quad \quad \text{and} \quad \quad
[\cmbar_{\GG}]^{\glued} =
[\cmbar_{\GG}]^{\vir}.$$}

\bpf We relate the obstruction theory for $f$ to the obstruction
theories for $f'$ and $f''$ separately.  We use the description of the
obstruction theory given in Section~\ref{pot}.  A priori, these
virtual classes are defined relative to two different bases.  The
first is defined relative to $\cT \times \fM_{g,m+n}$ and the
second is defined relative to $\cT \times \cT_\unrigid \times
\fM_{g',m'+\de} \times \fM_{g'',m''+\de + n}$.  However, the
torus action on the relative cotangent complex of the map between
these two base spaces has no torus fixed part, because the
$\C^*$-action on the deformation space of $D_0$ is nontrivial, as is
the torus action on the deformation spaces of the $\de$ nodes $\{ N_i
\}$.  Therefore we need only consider the relative obstruction
theories.  By \eqref{postit} and \eqref{localglobal} there are two
pieces here: one is $H^\cdot (C, f^\dagger T_X(-\log D))$, and the other is
the local obstruction at the nodes mapping to the singular locus of
$X_l$.  The local obstruction coming from the node $N_i$ has a
nontrivial $\C^*$-action, because the weight of this torus action is
the same as the weight of the action on $N_{D/X}$.  Thus, they do not
contribute to the fixed obstruction theory.  The local obstructions
coming from the singular locus of $Y_l$ have zero weight, since the
torus doesn't act on $Y_l$, meaning that these local obstructions
occur in the torus-fixed part of the perfect obstruction theory, just
as they do in the obstruction theory on $\M_{\Gamma''}$.  There are no
local obstructions for $f'$, since the target is smooth.  We study the
global obstructions using the partial normalization map $C' \coprod
C'' \rightarrow C$.  This gives us a long exact sequence in cohomology
$$ \xymatrix{0  \ar[r] &  H^0(C, f^\dagger T_{X_l}(-\log D_{\infty}))  } $$
$$ \xymatrix{ \ar[r] &  H^0(C', f'^* (T_X(-\log D))) \oplus
H^0(C'',f''^\dagger (T_{Y_l}(-\log D_0 - \log D_{\infty}))) } $$
$$ \xymatrix{ \ar[r] &  \oplus_{i=1}^n \left. T_D \right|_{f(N_i)} \ar[r] &
H^1(C, f^\dagger(T_{X_l}(-\log D_\infty)))  } $$
$$ \xymatrix{  \ar[r] &  H^1(C', f'^*(T_X(-\log D))) \oplus 
H^1(C'', f^\dagger (T_{Y_l}(-\log D_0 - \log D_{\infty})))  \ar[r] &  0.
 } $$ 
Note that this looks slightly different from the standard
normalization sequences, because the agreement required at the nodes
is only in the space $T_D$.  The torus action
on $T_D$ is trivial, since the divisor is torus-fixed.  We conclude that
the difference between the zero-weight piece of this obstruction
theory and the zero-weight piece of the obstruction theory coming from
$f'$ and $f''$ separately is precisely the term $\oplus \left. T_D \right|_{f(N_i)}$.
This is simply the pullback of the normal bundle of the diagonal
morphism from $D^n$ to $D^{2n}$ (which corresponds to the refined 
Gysin homomorphism $\De^{!}$ in \eqref{butter}).  Thus, this sequence is
exactly the verification of the compatibility condition of perfect
obstruction theories that guarantees that the two possible choices of
virtual class on this moduli space agree
 (see for example
\cite[Prop.\ 5.10]{bf}).  (This is essentially the
same as the 
standard argument used to prove the splitting axiom in Gromov-Witten
theory.)  \epf

\epoint{Two virtual bundles on $\cmbar_{\GG}$ and
$\cFbarGG$} 
\label{wine}\lremind{wine}We 
next describe two relevant (virtual) bundles on
$\cmbar_{\GG}$ and
$\cFbarGG$.  First, the pullback of $N_{\Gamma'}$ to
$\M_{\Gamma',\Gamma''}$ is a virtual bundle which we will also denote
by $N_{\Gamma'}$.  This bundle naturally 
descends to $\cFbarGG$.  For convenience, we denote the
descended bundle $N_{\Gamma'}$ as well.

The second relevant bundle is the line bundle $\cLL$ corresponding to
the deformation of the singularity $D_0$ ($=X \cap Y_l$ in $X_l$).
(This line bundle $\cLL$ is analogous to the line bundle $\mathbf{1}$
arising in J.~Li's degeneration formula \cite[p.~203]{li2}.) 
The fiber of this
bundle at a point of the boundary is canonically isomorphic to
$H^0(D,N_{D/X} \tensor N_{D_0/Y_l})$.  The line bundle $N_{D/X}
\tensor N_{D_0/Y_l}$ is trivial on $D$, so its space of global
sections is one-dimensional, and we can canonically identify this
space of sections with the fiber of the line bundle at any point $pt$
of $D$.  Thus we see that we can write the bundle $\cLL$ as a tensor
product of bundles pulled back from the two factors separately.  The
one coming from $\M_{\Gamma'}(X,D)$ is trivial, since it is globally
identified with $H^0(pt,\left. N_{D/X} \right|_{pt})$, but it has a
nontrivial torus action; we denote this weight (i.e.\ the first Chern
class of this bundle, which is pure weight) by $w$.  In other words,
$w$ is the weight of the torus action on the normal bundle to $D$ in
$X$.  The line bundle coming from $\M_{\Gamma''}(Y)$ is a nontrivial
line bundle, but with trivial torus action; this is precisely the
pullback of $-\psi$ (where $\psi$ was defined in Section~\ref{psi}).  Thus
$$
c_1(\cLL) = w - \psi.$$

\epoint{Remark} \label{flake}\lremind{flake}For each 
node $N_i$ joining $C'$
to $C''$,  
there is a natural isomorphism between the line bundle $\cLL$ and the $m_i^{\text{th}}$
tensor power of the line bundle corresponding to the deformation of the
singularity of $C$ at $N_i$.  This is  because the morphism from $C'$ 
(respectively $C''$) induces an isomorphism between 
$\left. T_{C'}^{\tensor m_i} \right|_{N_i}$
and $\left. N_{D/X} \right|_{f(N_i)}$ (respectively $\left. T_{C''}^{\tensor m_i} \right|_{N_i}$
and $\left. N_{D_0/Y} \right|_{f(N_i)}$).

\point We may now state the relative virtual localization theorem,
which reduces understanding the contributions of an
arbitrary fixed locus to understanding the contributions from a simple
fixed locus and the contributions from maps to a non-rigid target.


\tpoint{Theorem (Relative virtual localization)} \lremind{rvlt}
\label{rvlt}
{\em 
\begin{eqnarray*}
{} [\cmbar_{\Gamma}(X,D)]^{\vir} &=&  
 \frac { [\cmbar_{\Gamma}(X,D)^{\simple}]^{\vir}} { e(N_{\Gamma})  }
+ \sum_{\text{$\cmbar_{\GG}$ composite}}
\left(  \prod m_i \right) 
\frac  { [\cFbarGG]^{\vir}} { e(N_{\Gamma'}) c_1(\cLL)}
\\
 &=&   
 \frac { [\cmbar_{\Gamma}(X,D)^{\simple}]^{\vir}} { e(N_{\Gamma})  }
+ \sum_{\text{$\cmbar_{\GG}$ composite}}
\frac {\left(  \prod m_i \right) 
gl_* [\cmbar_{\Gamma',
  \Gamma''}]^{\glued}} { \left| \Aut( m_{\cdot} ) 
\right| e(N_{\Gamma'}) (w - \psi)}
\end{eqnarray*}
}

The two versions of the theorem are obviously equivalent by Lemma \ref{butterfly}.
Relative virtual localization may be interpreted (and proven) as follows.  First, 
the induced virtual fundamental
class on the fixed locus agrees with the natural virtual fundamental
class coming from the modular interpretation of the fixed locus
(Lemma~\ref{butterfly}).  Second,
there is a
contribution to the virtual normal bundle not present on the
moduli space of stable maps which arises from the deformation of
the target in the case of composite fixed loci.  This contribution is
virtual codimension $1$ (in a more precise sense virtual codimension
$(\de+1)-\de$, as we shall see in the proof), and the contribution to the
Euler class of the virtual normal bundle is $c_1(\cLL) / \prod m_i$.  

\epoint{Corollary: Relative virtual localization with target
  $\proj^1$} \label{bert}\lremind{bert}To prove \thmstar, we will use
  relative virtual localization for target $\proj^1$.  In this case we
  geometrically interpret $N_{\Gamma'}$ by giving the changes required
  from the formula of \cite[p.~505]{vl}.  We choose the $\C^*$-action
  on $\proj^1$ that acts with weight $1$ on $T_{\proj^1}$ at $0$ and
  weight $-1$ at $\infty$; for example, we can take $\la_0=t$ and
  $\la_{\infty}=0$ in the language of \cite{vl} (where $t$ is the
  generator of equivariant Chow ring of a point).  Thus $w=-t$.

Then the formula for $e(N_{\Gamma'})$ is the same as the formula
for $e(N^{\vir})$ of \cite[p.~505]{vl} except 
$\la_0 = t$ and $\la_{\infty} = 0$, and 
the edge contribution should be
$$
\prod_{\text{edges $e$}} \frac { d_e^{d_e}} { d_e! t^{d_e}}.
$$
This accounts for replacing $T_{\proj^1}$ with $T_{\proj^1}(-\infty)$.
There is no  correction to the vertex contributions: 
because there are no nodes
or contracted components over infinity, the other terms that could
conceivably change are unaffected by this replacement of bundles.

We now give two proofs of the Relative virtual localization
theorem~\ref{rvlt}.  The first is short and direct, using the
description of the perfect obstruction theory given in
Section~\ref{pot}.  The second proof avoids the technical details of
the perfect obstruction theory 
by invoking more results from \cite{li2} and independently establishes
the equality of virtual classes of Lemma \ref{butterfly}.  (This is possible because
the technical machinery required for relative virtual localization
is very similar to that needed for Li's degeneration formula.)

\bpoint{The first proof  of relative virtual localization}
To prove our relative virtual localization statement, what remains is
to analyze the virtual normal bundle, which is determined from the
non-fixed part of the perfect obstruction theory.  Part of this
virtual bundle is just the virtual normal bundle of the locus
$\M_{\Gamma'}$ inside $\M_{\Gamma'}(X,D)$.  (The $\C^*$-action on
$\M_{\Gamma''}(Y,D_0, D_{\infty})_\unrigid$ is trivial, so there is no
virtual normal bundle on that side.)  All that remains are the terms
coming from the nodes $\{ N_i \}$ of $C$ and the singularities of $X_l$.  There is a
single deformation of $X_l$ with nonzero torus weight corresponding to
smoothing the singularity $D_0$ of $X_l$.  This contributes a copy
of $\cLL$ to the virtual normal bundle.  For each of the $N_i$ there
is a deformation of $C$ corresponding to smoothing the node at $N_i$,
giving a contribution of $\frac{1}{m_i} \cLL$ (Remark~\ref{flake}).  Finally, each
$N_i$ also contributes a local obstruction which we have seen is given
by $\cLL$ (see \eqref{explocalglobal}).  We cancel the resulting copy of $c_1(\cLL)$ that
occurs in both the numerator and denominator $n$ times, leaving
$\frac{1}{\prod m_i}$.  We conclude that the difference between the
Euler class of the virtual normal bundle to $\M_{\Gamma',\Gamma''}$ in
$\M_{\Gamma}(X,D)$ and the Euler class of the pullback of the virtual
normal bundle of $\M_{\Gamma'}$ in $\M_{\Gamma'}(X,D)$ is given 
by $\frac{1}{\prod m_i}c_1(\cLL)$, which completes the proof of Theorem
\ref{rvlt}. \epf

\bpoint{The second proof of relative virtual localization}
We again show that the Euler class of the virtual normal bundle to $\cFbarGG$
is $e(N_{\Gamma'})
(w - \psi)/ \prod m_i$.
We first define two virtual divisors on $\cmbar_{\Gamma}(X,D)$, 
denoted $Z_{\thick}$ and $Z_{\thin}$.  

There is a line bundle $\cLL'$ with section $s_{\cLL'}$ on
$\cmbar_{\Gamma}(X,D)$ corresponding to the locus where the target
breaks, and the source curve data splits into $\Gamma'$ and
$\Gamma''$.  Let $Z_{\thick}$ be the virtual Cartier divisor
corresponding to $s_{\cLL'}$.  \'{E}tale-locally near a point of
$Z_{\thick}$, the divisor is pulled back from the deformation space of
the singularity $D_{\GG}$ of $X_l$ splitting the
source curve data into $\Gamma'$ and $\Gamma''$.  Then $Z_{\thick}$
has an induced obstruction theory, and hence a virtual fundamental
class given by
$$
[Z_{\thick}]^{\vir} = c_1(\cLL') \cap [\cmbar_{\Gamma}(X, D)]^\vir$$
by
\cite[Lemma~3.11]{li2} (see also \cite[Lemma~4.6]{li2}).

The restriction of $\cLL'$ to $\cFbarGG$ is the line bundle $\cLL$ of Section~\ref{wine}.  ({\em
  Caution:} Although $c_1(\cLL) = w- \psi$, $c_1(\cLL') \neq w - \psi$
on $\cmbar_{\Gamma}(X,D)$ in general, because $Z_\thick$ may
contain maps corresponding to points $\cmbar_{\Gamma'}(X,D_{\GG})$ where the
target $X$ degenerates.  Then
the normal bundle $N$ to the $D_{\GG}$ in the degeneration of $X$ is
not the same equivariantly as $N_{D/X}$.)

Define $Z_{\thin}$ as the closed substack of $Z_{\thick}$ where we
require furthermore that the nodes of the source curve mapping to
$D_{\GG}$ are not smoothed (even to first order):  
$$
Z_{\thin} = gl \left( \cmbar_{\Gamma'}(X,D) \times_{D^n}
\cmbar_{\Gamma''}(Y,D_0, D_{\infty})_{\unrigid} \right),
$$
where $gl$ again is the \'{e}tale quotient map by the finite
group $\Aut ( m_\cdot)$.  $Z_{\thin}$ is given an obstruction theory
in \cite[p.~250--252]{li2}.  Let $[Z_{\thin}]^{\vir}$ be the
associated virtual fundamental class.  A second virtual fundamental
class (and obstruction theory) comes from the gluing description:
$$
[Z_{\thin}]^{\glued} =  
\frac 1 {\left| \Aut (m_\cdot) \right|} gl_* \Delta^!
\left(
[\cmbar_{\Gamma'}(X,D)]^{\vir} \times_{D^n}
[\cmbar_{\Gamma''}(Y,D_0, D_{\infty})_{\unrigid}]^{\vir} \right)
$$
(compare to \eqref{butter}).  By \cite[Lemma~3.14]{li2},
$[Z_{\thin}]^{\glued} = [Z_{\thin}]^{\vir}$, and in fact the two
obstruction theories are identical \cite[Lemma~4.15]{li2}.

Consider the sequence of inclusions
\begin{equation} \label{candle}
\cFbarGG 
\subset
 Z_{\thin}
\subset
Z_{\thick}
\subset
\cmbar_{\Gamma}(X,D).
\end{equation}
(The last three terms should be compared to the first, second,
and fourth terms of the sequence of inclusions of \cite[p.~248]{li2}.)
For each inclusion, we identify the virtual normal bundle of each 
term in the next, and verify that it has no fixed part when restricted
to $\cFbarGG$.
As the obstruction theory of $\cFbarGG$ is
$\C^*$-fixed, 
it follows that
the obstruction
theory of $\cFbarGG$ is indeed the fixed part of
the restriction of the obstruction theory of $\cmbar_{\Gamma}(X,D)$.
By multiplying the Euler classes of the virtual normal bundles of each inclusion,
we obtain the Euler class of the virtual normal bundle of 
$\cFbarGG$, completing the proof of relative
virtual localization.

\noindent
{\em The inclusion $Z_{\thick} \subset \cmbar_{\Gamma}(X,D)$.}  As
stated earlier, the virtual normal bundle of $Z_{\thick}$ in
$\cmbar_{\Gamma}(X,D)$ is $\cLL'$.  When restricted to
$\cFbarGG$, it is $\cLL$, which has no fixed part,
as $\C^*$ acts on $\cLL$ nontrivially (by our assumption that $\C^*$
does not act trivially on $X$).

\noindent
{\em The inclusion $Z_{\thin} \subset Z_{\thick}$.}  We need to delve
into Li's argument of \cite[Sect.~4.4]{li2}, where the obstruction
theories of $Z_{\thin}$ and $Z_{\thick}$ are compared, to verify that
the fixed parts of both obstruction theories (when restricted to
$\cFbarGG$) are the same.  We refer the reader in
particular to the first paragraph of that section for an overview of
the strategy.

We motivate the argument with a simple example.  Consider the
morphism of schemes from the formal neighborhood of $\de$ nodes
to the formal neighborhood of a node, where
the morphism is predeformable, with the $i^{\text{th}}$ node of the source
mapping with ramification $m_i$.  The deformation space of this
morphism is reducible if $\prod m_i >1$: if (to first order) any of
the branchings of the source curve move away from the
(scheme-theoretic) preimage of the target node, then the target node
cannot smooth, even to first order.  There is precisely one component
that surjects onto the deformation space of the target node.  Formal
equations for this component are\lremind{shaker}
\begin{equation} \label{shaker}
\xymatrix{
\Spf \C[[y_1, \dots, y_\de, x]] / (x = y_1^{m_1} = \cdots = y_\de^{m_\de})
\ar[r] &  \Spf \C[[x]].}
\end{equation}
Here $x$ corresponds to the deformation parameter of
the target node, and $y_i$ to the deformation parameter of the
$i^{\text{th}}$ node of the source curve.
(See \cite[Sect.~2.5]{germ} for this deformation-theoretic calculation.
This argument can be extended to higher-dimensional targets, 
and this has been done by many authors; 
see for example \cite{ch,v1} in the algebraic category.)

A consequence of J. Li's obstruction theory is that the relationship between
$Z_{\thick}$ and $Z_{\thin}$ is ``virtually'' analogous to
\eqref{shaker}.
More precisely, consider the formal thickening of $Z_{\thin}$ by the
formal deformation space of $D_{\GG}$ and the formal
deformation spaces of the nodes $\{ N_i \}$ of the source curve mapping to the
singularity of $D_{\GG}$.  Locally on $Z_{\thin}$, $x$
is a generator of the first deformation space, and $y_1$, \dots, $y_\de$
are generators of the rest.  Then $Z_{\thick}$ is the pullback of
$x=y_1^{m_1} = \cdots = y_\de^{m_\de}=0$, and $Z_{\thin}$ is the pullback
of $x=y_1 = \cdots = y_\de=0$.  Li describes these line bundles
explicitly, but this description is unnecessary for our argument
(even though the torus acts nontrivially on them).  Both $Z_{\thick}$
and $Z_{\thin}$ sit in a space of virtual dimension $\de$ larger,
corresponding (again \'etale-locally) to
$$
\Spec \oh_{Z_{\thin}} \hookrightarrow \Spec \oh_{Z_{\thin}}[y_1,
\dots, y_\de, x] / (x = y_1^{m_1} = \cdots = y_\de^{m_\de}) \hookrightarrow
\Spec \oh_{Z_{\thin}}[[y_1, \dots, y_\de, x]].$$ The larger space has a
natural obstruction theory. The obstruction theory of $Z_{\thin}$ is
obtained by capping with $\de$ Cartier divisors $y_1 = \cdots = y_\de =
0$, and the obstruction theory of $Z_{\thick}$ is obtained by capping
with $\de$ Cartier divisors $y_1^{m_1} = \cdots = y_\de^{m_\de} = 0$.  (This
requires elaboration, as one cannot do this for a general perfect
obstruction theory.  For a justification, see \cite[Sect.~4.4]{li2}.)

Thus the fixed part of the obstruction theory of $Z_{\thin}$
(restricted to $\cFbarGG$) agrees with the fixed
part of the obstruction theory of the larger space (restricted to
$\cFbarGG$), which agrees with the fixed part of the
obstruction theory of $Z_{\thick}$ (restricted to $\cFbarGG$).

This analysis (or \cite[Lemma~3.12]{li2}) also implies that
$[Z_{\thin}]^{\vir} = \frac 1 {\prod m_i} [Z_{\thick}]^{\vir}$, so the
relative codimension of $Z_{\thin}$ in $Z_{\thick}$ is $0$, and the
Euler class of the virtual normal bundle is $1 / \prod m_i$.

We remark that $Z_{\thick}$ is thus virtually Cartier
on $\cmbar_{\Gamma}(X,D)$, locally cut out by one equation, but that
$Z_{\thin}$ is not: in some sense it is cut out by $(\de+1)-\de$
equations.

\noindent
{\em The inclusion $\cFbarGG \subset Z_{\thin}$.}
We use the gluing version of the obstruction theory of $Z_{\thin}$.
By pullback of obstruction theories, the virtual normal bundle
of
$$\cmbar_{\Gamma'}(X,D)^{\simple} \times_{D^n} \cmbar_{\Gamma''}(Y, D_0, D_{\infty})_{\unrigid} \quad \text{ in } \quad
\cmbar_{\Gamma'}(X,D) \times_{D^n} \cmbar_{\Gamma''}(Y, D_0, D_{\infty})_{\unrigid}$$
is the pullback of the virtual normal bundle of
$\cmbar_{\Gamma'}(X,D)^{\simple}$ in
$\cmbar_{\Gamma'}(X,D)$, i.e.\ $N_{\Gamma'}$.
By descending by $gl$, we see that the normal bundle of
$\cFbarGG$ in $Z_{\thin}$ is $N_{\Gamma'}$.
By definition of the virtual normal bundle of $\cmbar_{\Gamma}(X,D)^{\simple}$
in $\cmbar_{\Gamma'}(X,D)$ (i.e.\ arising from the non-fixed part of the
obstruction theory), the virtual normal bundle
has no fixed part.

In conclusion, the obstruction theory of $\cFbarGG$
is indeed the fixed part of the restriction of the obstruction theory 
of $\cmbar_{X,D}$, and the Euler class of its virtual normal bundle
is $(w- \psi) \times (1 / \prod m_i) \times e(N_{\Gamma'})$ as desired,
completing the proof of relative virtual localization (Theorem~\ref{rvlt}).
\epf

\section{Vanishing of tautological classes on moduli
spaces of curves and stratification by number of rational components
(\thmstar)}
\label{starsection}
\lremind{starsection}
In this section, we give background on the moduli space of curves
and its tautological ring.  We then define {\em Hurwitz cycles}
(in the Chow ring of the moduli space of curves), and show that 
Hurwitz cycles satisfy the conclusion of \thmstar, using
degeneration techniques.  
The proof of \thmstar {} will then follow by showing that
tautological classes are essentially linear combinations
of Hurwitz cycles, using relative virtual localization.

\bpoint{Background on the moduli space of curves} (See \cite{notices}
for a more leisurely survey of the facts we will need about the moduli
space of curves and its tautological ring.)  We assume familiarity
with $\mgn$.  

There are {\em natural morphisms} among moduli spaces
of stable curves,
forgetful morphisms
\begin{equation}
\label{forgetful}
\xymatrix{ \cmbar_{g,n} \ar[r] &  \cmbar_{g,n-1}}
\end{equation}
and gluing morphisms\lremind{gluing1, gluing2}
\begin{eqnarray}   \label{gluing1}
 \cmbar_{g_1, n_1+1} \times 
\cmbar_{g_2, n_2 +1}  & \xymatrix{  \ar[r] &  } &  \cmbar_{g_1+g_2, n_1+n_2},   \\
\label{gluing2}
\cmbar_{g, n +2} 
  & \xymatrix{ \ar[r] & } &  \cmbar_{g+1, n}. 
\end{eqnarray}
Gluing morphisms will be denoted by $gl$.
The tautological ring may be defined in terms of the natural morphisms
as follows.

\epoint{Definition}  \lremind{tautdef}\label{tautdef}The system of {\em tautological rings} are defined as the smallest
system of $\Q$-vector spaces $\left( R^i \left( \cmbar_{g,n} \right) \right)_{i,g,n}$ satisfying:
\begin{itemize}
\item $\psi_1^{a_1} \cdots \psi_n^{a_n} \subset R^*\left( \mgn \right)$, and
\item the system is closed under pushforwards by the natural morphisms.
\end{itemize}
(Throughout we consider Chow groups with $\Q$-coefficients.)
The tautological ring of a dense open subset of $\mgn$ is defined by
restriction.  Let $R_j\left( \mgn \right)= R^{\dim \mgn-j}\left( \mgn \right)$
be the group of dimension $j$ tautological classes.

This is equivalent to the other definitions of the tautological ring
appearing in the literature.  For example, to show equivalence
with the definition of \cite{nontaut}, by
\cite[Prop.~11]{nontaut} it suffices to show that any 
monomial in the $\psi$-classes and $\kappa$-classes lies
in the groups defined in Definition~\ref{tautdef}.  
But such
a class is clearly the pushforward of a monomial in $\psi$-classes
on a space of curves with more points (as in \cite[p.~413]{lthm}).

Let $\cS_k$ (resp.\ $\cS_{\geq k}$, etc.) be the union of strata of
$\cmbar_{g,n}$ with precisely $k$ (resp.\ at least $k$, etc.)
components of geometric genus $0$.  Then $A_*(\cS_{\geq k}) \rightarrow
A_*\left( \mgn \right) \rightarrow A_*(\cS_{<k}) \rightarrow 0$ is exact.  Let $I_{
  \geq k}$ be the image of $A_*(\cS_{\geq k})$ in $A_*\left( \mgn \right)$; it is
an ideal of $A^*\left( \mgn \right)$.  Then \thmstar{} is equivalent to
$R^i\left( \mgn \right) \subset I_{ \geq i-g+1}$.  (It will be useful 
to observe that the image of $I_{ \geq k}$ under a forgetful morphism
is $I_{ \geq k-1}$.)

Define $\cmt_{g,n}$ in the same way as $\cmbar_{g,n}$, except the
curve is not required to be connected.  Then $\cmt_{g,n}$ is a
Deligne-Mumford stack of finite type, and is stratified by $\cS_k$
analogously; each irreducible component is isomorphic to the quotient
of a product of $\cmbar_{g',n'}$'s by a finite group.  Let $\cmt =
\coprod_{g,n: 2g-2+n > 0} \cmt_{g,n}$.

\epoint{Remark} \label{ernie}\lremind{ernie}Note 
that \thmstar {} is true for given $g,n$ with
$\cmbar$ replaced by $\cmt$ if it is true for all $\cmbar_{g',n'}$
with $\dim \cmbar_{g',n'} \leq \dim \cmbar_{g,n}$.  Hence \thmstar {}
implies that the corresponding statement with $\cmbar$ replaced by
$\cmt$ is true in general.

\bpoint{Hurwitz classes, and their behavior with respect to the
  stratification} \label{bub}\lremind{bub}For the rest of
Section~\ref{starsection}, fix $g$ and $n$.  Following Ionel, define a
{\em trivial cover} of $\left( \proj^1, 0, \infty \right)$ to be a map
of the form $\proj^1 \rightarrow \proj^1$, $[x;y] \mapsto [x^u;y^u]$,
or equivalently a map from an irreducible curve to $\left( \proj^1, 0,
  \infty \right)$ with no branch point away from $0$ and $\infty$.  A
{\em trivial component} of a stable relative map is a connected
component of the source curve, such that the stabilization of the
morphism from that component is a trivial cover.

Let $\cmbar_{g,\al}\left( \proj^1 \right)$ be the stack
parametrizing degree $d$ stable relative maps from a curve 
of arithmetic genus $g$ to $\proj^1$, relative
to one point $\infty$, corresponding to partition $\al \vdash d$.
(As stated earlier, curves are not assumed to be connected, and
partitions are taken to be ordered, i.e.\ the parts are
labeled.)  Let $\cmbar_{g,\al, \be}\left( \proj^1 \right)$ be the moduli space of
stable relative maps relative to two points $0$ and $\infty$ (with
corresponding partitions $\al \vdash d$ and $\be \vdash d$
respectively).  Let $\cmbar_{g,\al,\be}\left( \proj^1 \right)_{\unrigid}$ be the
analogous space of maps to a non-rigid $\proj^1$
(Section~\ref{nonrigid}).

We will make repeated use of the following result.

\tpoint{Theorem (Degeneration formula, special case of \cite[Thm.\ 3.15]{li2})}\label{degthm}
{\em 
$$
[ \cmbar_{g,\al,\be}(\proj^1)]^{\vir} = \bigoplus \frac { \prod m_i} {\Aut(m_\cdot)}
gl_* \left(  [ \cmbar_{g', \al, m_\cdot}(\proj^1) ]^{\vir} \boxtimes
[ \cmbar_{g'', m_\cdot, \be}(\proj^1) ]^{\vir} \right),$$
where the sum on the right is over all possible choices of $g'$, $g''$,
$m_\cdot$, and the $\al$ and $\be$ may be omitted.
}

This equality should be interpreted in the total space of 
the family of relative stable map spaces associated to a family of
$\proj^1$'s degenerating to a pair of $\proj^1$'s meeting in a point.
In practice, we will use it by capping these virtual classes against 
natural Chow cohomology classes that make sense for such a family, and 
pushing the resulting equality forward to a moduli space of curves.

(We are grateful to Y.~Ruan for pointing out that in the analytic setting,
one can read off the degeneration formula for virtual fundamental
classes from the first four lines of the (earlier) proof of Theorem 5.6
in \cite{liruan}.  He further remarks that the virtual neighborhood is
a much stronger notion than the virtual fundamental class because
it is a `smooth object'.)

Let $r = \vdim \left( \cmbar_{g,\al}\left( \proj^1 \right) \right) = d+n+2g-2$ (where $\vdim$
denotes {virtual dimension}). Note that $r$ is the dimension of the
open set $\cm_{g,\al}\left( \proj^1 \right)$ parametrizing maps with smooth source
(the Hurwitz scheme).  There is a
branch morphism $br:
\cmbar_{g,\al}\left( \proj^1 \right) \rightarrow \Sym^r \proj^1$, obtained by
considering the morphisms in $\cmbar_{g,\al}\left( \proj^1 \right)$ as morphisms to
$\proj^1$, and noting that the image of the Fantechi-Pandharipande
branch morphism \cite[Thm.~1]{fnp} to $\Sym^{r+\sum (\al_i-1)}
\proj^1$ always contains $\infty$ with multiplicity at least $\sum (\al_i-1)$.

Define a pushforward morphism $\pi$ from $A_* \left(
  \cmbar_{g,\al,\be}\left( \proj^1 \right) \right)$ (and $A_* \left(
  \cmbar_{g,\al,\be}\left( \proj^1 \right)_{\unrigid}\right)$) to
$A_*\left( \cmt \right)$ as follows.  On the locus of
$\cmbar_{g,\al,\be}$ where the cover has no trivial parts, $\pi$ is
the usual pushforward (by the morphism induced by the universal curve;
the marked points correspond to the parts of $\al$ and $\be$).  On the
locus where the cover has trivial parts, ignore the points on the
trivial parts and push forward.  (Note that the automorphism groups of
the trivial covers still play a role, as they contribute a multiplicity
to the virtual fundamental class $[\cmbar_{\al, \be}(\proj^1)]^\vir$: the
inverse of the product of their degrees.)  This definition is designed
so that the degeneration formula for fundamental classes
(Theorem~\ref{degthm}) agrees with the gluing of strata (used in
Lemma~\ref{l3} and Proposition~\ref{firstprop} below).

Now assume $r \geq j$.  Define the {\em Hurwitz class} $\HH^{g,\al}_j$
informally by considering $[\cmbar_{g,\al}\left( \proj^1 \right)]^{\vir}$, fixing
$r-j$ branch points (leaving a class of dimension $j$), and pushing
forward to $\cmt$.  More precisely, $$\HH^{g,\al}_j = \pi_* \left(
\cap_{s=1}^{r-j} br^*(L_{p_s}) \cap [\cmbar_{g,\al}\left( \proj^1 \right)]^{\vir}
 \right),$$
where $p_1, \dots, p_{r-j} \in \proj^1$, and $L_p$ corresponds to the
hyperplane in $\Sym^r\left( \proj^1 \right)$ of $r$-tuples containing $p$.

Define the {\em double Hurwitz class} $\HH^{g,\al,\be}$ by considering
$[\cmbar_{g,\al,\be}\left( \proj^1 \right)]^{\vir}$, 
fixing one branch point, and pushing forward to $\cmt$:
$$\HH^{g,\al,\be} = \pi_* \left( 
br^*(L_1) \cap [\cmbar_{g,\al, \be}\left( \proj^1 \right)]^{\vir} \right)\in A_*\left( \cmt \right).$$
  Define
$\HH^{g,\al,\be}_{\unrigid}$ similarly by 
$$
\HH^{g,\al,\be}_{\unrigid} = \pi_*
\left[ \cmbar_{g,\al,\be}\left( \proj^1 \right)_{\unrigid} \right]^{\vir} \in A_*\left( \cmt \right).$$  
Then
$\dim \HH^{g,\al,\be} = \dim \HH^{g,\al,\be}_{\unrigid} =
r(g,\al,\be)-1$, where $r(g,\al,\be) =\vdim \left( \cmbar_{g,\al,\be}\left( \proj^1 \right) \right)$
is the ``expected number'' of branch points away from $0$ and $\infty$.

Note that $\HH^{g, \al}_j$ is the
pushforward of a class on $\cap_{s=1}^{r-j} br^{-1}(L_{p_s}) \cap
\cmbar_{g,\al}\left( \proj^1 \right)$, and similarly for the double Hurwitz class $\HH^{g,\al,\be}$.

\tpoint{Lemma} {\em $\HH^{g,\al,\be} = r(g,\al,\be) \HH^{g,\al,\be}_{\unrigid}$.}
\lremind{comparison} \label{comparison}

\bpf (The commutative diagram below may be helpful.)  Let $T$ be the
smooth locus of the universal target over $\cmbar_{g,\al,\be}\left(
\proj^1 \right)_{\unrigid}$, minus the $0$- and $\infty$-sections.
Let $B$ the branch locus of the universal map (a subset of $T$).  Then
$B$ is an effective Cartier divisor on $T$, finite of degree $r(g,\al,\be)$
over the base (the arguments of \cite[Sect.~3]{fnp} for stable maps
apply without change for stable relative maps).  Then $T$ is
canonically an open subset of $\cmbar_{g,\al,\be}\left( \proj^1
\right)$ (using the three points $0$, $\infty$, and the point on the
universal target to give a morphism to $\proj^1$, sending the
universal point to $1$), and their obstruction theories are identical.
Indeed, this open set of $\M_{g,\alpha,\beta}\left( \proj^1 \right)$
is just the corresponding open subset of the base-change of
$\M_{g,\alpha,\beta}\left( \proj^1 \right)_\unrigid$ by the forgetful
morphism from $\cT$ to $\cT_\unrigid$.  Since the virtual classes are
constructed from relative perfect obstruction theories over $\cT$ and
$\cT_\unrigid$ respectively, this follows from the base change
property of the virtual class construction. Thus, the restriction of
$[\M_{g,\alpha,\beta}\left( \proj^1 \right)]^\vir$ to $T$ is just the
flat pull-back of $[\M_{g,\alpha,\beta}\left( \proj^1
\right)_\unrigid]^\vir$.
Under this open immersion, the effective
Cartier
divisor $B$ on $T$ is canonically the effective Cartier divisor
$br^{-1}(L_1)$.   
Since the degree of the morphism $B\rightarrow \M_{g,\alpha,\beta}\left( \proj^1 \right)_\unrigid$
is $r(g,\alpha,\beta)$, the result follows.
$$
\xymatrix{
\cmbar_{g,\al,\be}\left( \proj^1 \right)  & &  br^{-1}(L_1) 
\ar@{_{(}->}[ll]_{\txt{\scriptsize{Cartier}}} 
\\
T \ar[u]^{\txt{\scriptsize{open imm.}}} 
\ar[dr] & & B   \ar[dl]^{\txt{\scriptsize{ \quad \quad finite, degree $r(g,\al,\be)$}}}
\ar@2{-}[u]
\ar@{_{(}->}[ll]_{\txt{\scriptsize{Cartier}}}  \\
&   \cmbar_{g,\al,\be}\left( \proj^1 \right)_{\unrigid}
} 
$$
\epf

\tpoint{Lemma} {\em If $(g,\al,\be) \neq (0, (d), (d))$, then
$\HH^{g,\al,\be} \in I_{\geq 1}$.} \label{ilem}
\lremind{ilem}

We are grateful to E.-N. Ionel for pointing out this result to us.
The proof below is extracted from her proof of
\cite[Prop.~2.8]{igetz}, with the following differences: this argument
is in the virtual setting; the source curve is not required to be
connected, and the Ionel-Parker gluing formula is replaced by the J.
Li degeneration formula (Theorem~\ref{degthm}).  In particular, the
key trick in the following
argument (to use a surprising forgetful morphism)  is due to her.

\bpf We argue by induction on $\dim \HH^{g,\al,\be}$, and then on the
degree $d= \left| \al \right|$.  Assume that we know the result for
$(g',\al',\be')$ with $\dim \HH^{g',\al',\be'} < \dim
\HH^{g,\al,\be}$, and for $\dim \HH^{g',\al',\be'}=\dim
\HH^{g,\al,\be}$ but $\left| \al' \right|< \left| \al \right|$.
Consider $ev^{-1}(p) \cap \cmbar_{g,\al,\be,1}(\proj^1)$,
parametrizing relative stable maps of the sort we are interested in,
with a branch point fixed at $1$, and with an additional marked point
$q'_1$ mapping to some fixed  $p \in \proj^1$.  Let $\HH$ be its
virtual fundamental class.  Let $\rho$ be the forgetful morphism {\em
forgetting $q_1$, the marked point corresponding to $\al_1$}.  By
breaking the class $\HH$ in two ways, we will show that $\rho_* \pi_*
\HH \in I_{\geq 1}$, and that $\rho_* \pi_* \HH$ is a non-zero
multiple of $\HH^{g,\al,\be}$ modulo $I_{\geq 1}$, from which the
result follows.

Degenerate $\HH$ by breaking the target into two pieces, so that $0$
and $1$ lie on $T_1$, say, and $p$ and $\infty$ lie on $T_2$.  By the
degeneration formula (Theorem~\ref{degthm}), we can express $\HH$ as
the sum of other double Hurwitz classes, appropriately glued.
Consider one such summand $gl( \HH_1 \boxtimes \HH_2)$, where $\HH_i$
corresponds to the stable relative map to $T_i$.  The cover of $T_1$
has branching away from the two special points $0$ and $T_1 \cap T_2$,
as it has branching above $1$.  If the cover of $T_2$ also has
branching away from its two special points $T_1 \cap T_2$ and
$\infty$, then by the inductive hypothesis $\pi_*(\HH_1) \in I_{\geq
1}$, and also $\pi_*(\HH_2) \in I_{\geq 1}$ (by applying the inductive
hypothesis to the map without the marked point $q'_1$, and then
pulling back by the forgetful morphism; we use here the pullback property of the virtual
fundamental class).  Thus $\pi_*(\HH) \in I_{\geq 2}$, from which
$\rho_* \pi_* (\HH) \in I_{\geq 1}$.  On the other hand, if the cover
of $T_2$ has no branching away from the two special points, then we
can immediately identify a three-pointed rational curve in the
stabilization of the source: the component mapping to $T_2$ containing
$q'_1$, one of the marked points mapping to $\infty$, and a node
mapping to $T_1 \cap T_2$.  The one exception is if this node over
$T_1 \cap T_2$ is glued to a trivial cover of $T_1$, with the point
above $0$ marked $q_1$, as then once $q_1$ is forgotten, this
component is stabilized away.  But in this case,  the inductive
hypothesis for covers of lower degree implies that the contribution
from the remainder of the curve mapping to $T_1$ lies in $I_{\geq 1}$.
Hence we have shown that $\rho_* \pi_* \HH \in I_{\geq 1}$.

Next, degenerate $\HH$ by breaking the target so that $0$ and $p$ lie
in $T_1$ and $1$ and $\infty$ lie in $T_2$.  We proceed as in the
previous paragraph.  Consider any summand  $gl(\HH_1 \boxtimes \HH_2)$
in the degeneration
formula.  If there is branching over $T_1$
away from the two special points, then this summand lies in $I_{\geq
  1}$ by the same argument.  If there is no other branching over
$T_1$, and $q'_1$ lies in a component of the source containing another
marked point $q_i$ above $0$ {\em where $i \neq 1$}, then as before we
can identify a rational curve in the stabilization of the source, so
this summand lies in $I_{ \geq 1}$ again.  The only remaining case is
if there is no branching over $T_1$, and $q'_1$ lies in the component
of the source containing $q_1$; in this case, the contribution is
precisely $\HH^{g,\al,\be}$ times the degree of the trivial cover of
$T_1$ containing $q'_1$ and $q_1$. Hence $\rho_* \pi_* \HH$ is a
non-zero multiple of $\HH^{g,\al,\be}$ modulo $I_{\geq 1}$, from which
the result follows.  \epf

\tpoint{Lemma}  
{\em $\psi^k \cap \HH^{g,\al,\be}_{\unrigid} \in I_{\geq k+1}$.}
\label{l3}\lremind{l3}

Here $\psi$ is the class defined in Section~\ref{psi}, and $\psi^k
\cap \HH^{g,\al,\be}_{\unrigid}$ should be interpreted as $\pi_*
\left( \psi^k \cap \left[ \cmbar_{g,\al,\be}\left( \proj^1 \right)_{\unrigid}
  \right]^{\vir} \right)$.

\bpf We prove the result by induction.  The case $k=0$ follows from
the previous two lemmas.  Assume now that $k>0$.  By
Lemma~\ref{comparison} it suffices to show the analogous result on
$br^{-1}(L_1) \cap \cmbar_{g,\al,\be}\left( \proj^1 \right)$.  Via the target (with
three marked sections $0$, $1$, $\infty$), this space admits a
morphism to $\overline{\mathfrak M}_{0,3}$, and $\psi$ is pulled back from
$\psi_0$ on $\cT \subset \overline{\mathfrak{M}}_{0,3}$.  Now $\psi_0$ is
equivalent to a boundary divisor.
(This is {\em not} true on $\overline{\mathfrak
  M}_{0,2}$, hence the necessity of lifting to
$\cmbar_{g,\al,\be}\left( \proj^1 \right)$ instead of working on
$\cmbar_{g,\al,\be}\left( \proj^1 \right)_{\unrigid}$!)
More precisely, $\psi_0$ is the boundary divisor corresponding to (degenerations of) 
the
curves shown in Figure~\ref{orchid}, so
$$
\psi^k \cap \HH^{g,\al,\be}_{\unrigid} =
\sum gl \left(  \left( 
\psi^{k-1} \cap \HH^{g', \al', \be'}_\unrigid \right) \boxtimes \HH^{g'', \al'', \be''}
\right).
$$
But $\psi^{k-1} \cap \HH^{g', \al', \be'}_\unrigid \in I_{\geq k}$
by the inductive hypothesis, and $\HH^{g'', \al'', \be''} \in I_{\geq 1}$ by
Lemma~\ref{ilem}, so we are done, as the degeneration formula
(Theorem~\ref{degthm}) is compatible under $\pi$ with gluing.\epf

\begin{figure}
\begin{center}
\setlength{\unitlength}{0.00083333in}
\begingroup\makeatletter\ifx\SetFigFont\undefined%
\gdef\SetFigFont#1#2#3#4#5{%
  \reset@font\fontsize{#1}{#2pt}%
  \fontfamily{#3}\fontseries{#4}\fontshape{#5}%
  \selectfont}%
\fi\endgroup%
{\renewcommand{\dashlinestretch}{30}
\begin{picture}(1674,411)(0,-10)
\put(1212,162){\blacken\ellipse{50}{50}}
\put(1212,162){\ellipse{50}{50}}
\put(1437,87){\blacken\ellipse{50}{50}}
\put(1437,87){\ellipse{50}{50}}
\put(462,162){\blacken\ellipse{50}{50}}
\put(462,162){\ellipse{50}{50}}
\path(12,12)(912,312)
\path(762,312)(1662,12)
\put(1212,312){\makebox(0,0)[lb]{\smash{{{\SetFigFont{8}{9.6}{\rmdefault}{\mddefault}{\updefault}$1$}}}}}
\put(1512,237){\makebox(0,0)[lb]{\smash{{{\SetFigFont{8}{9.6}{\rmdefault}{\mddefault}{\updefault}$\infty$}}}}}
\put(312,312){\makebox(0,0)[lb]{\smash{{{\SetFigFont{8}{9.6}{\rmdefault}{\mddefault}{\updefault}$0$}}}}}
\end{picture}
}
\end{center}
\caption{$\psi_0$ on $\cT \subset \fM_{0,3}$\lremind{orchid}}
\label{orchid}
\end{figure}

We note that this proof can be unwound to show that $\psi^k \cap
\HH^{g, \al, \be}_{\unrigid}$ is equivalent to a sum of classes
corresponding to covers of a chain of $\proj^1$'s.

\tpoint{Proposition}  {\em $\HH^{g,\al}_j \in I_{\geq 2g-2+n-j}$.}
\label{firstprop}\lremind{firstprop}

\bpf Choose the points $\{ p_s \}$ to be distinct.
The class ${\mathbf
  H} = br^*( \cap_{s=1}^{r-j} L_{p_s})\cap
[\cmbar_{g,\al}\left( \proj^1 \right)]^{\vir}$ is supported on $br^{-1}( \cap
L_{p_s})\cap [\cmbar_{g,\al}\left( \proj^1 \right)]$.  Degenerate the target into a
chain of $r-j$ components $T_1$, \dots, $T_{r-j}$, with $p_s \in
T_s$,  and $\infty \in T_{r-j}$ (see Figure~\ref{oscar}). By the
degeneration formula (Theorem~\ref{degthm}), capped with
$br^*(L_{p_s})$, supported on $br^{-1}(L_{p_s})$), and the
compatibility of the degeneration formula with $\pi$, $\HH^{g,\al}_j =
\pi_* \mathbf{H}$ equals a sum with terms of the form\lremind{grover}
\begin{equation}\label{grover}
\pi_* \left( gl \left( 
  \boxtimes_{s=1}^{r-j} br^*(L_{p})[\cm_s]^{\vir}\right) \right)
=
 gl \left( 
  \boxtimes_{s=1}^{r-j} \pi_*\left( br^*(L_{p})[\cm_s]^{\vir}\right)\right) 
,
\end{equation}
where
$\cm_s$ is the appropriate moduli space of stable relative maps to $T_s$.
The product is taken over all appropriate matching of points.  Note that
while $\cm_1$ is the space of maps to $\proj^1$ relative to one point, for
all $s>1$, $\cm_s$ is a space of maps to $\proj^1$ relative to two points.

\begin{figure}
\begin{center}
\setlength{\unitlength}{0.00083333in}
\begingroup\makeatletter\ifx\SetFigFont\undefined%
\gdef\SetFigFont#1#2#3#4#5{%
  \reset@font\fontsize{#1}{#2pt}%
  \fontfamily{#3}\fontseries{#4}\fontshape{#5}%
  \selectfont}%
\fi\endgroup%
{\renewcommand{\dashlinestretch}{30}
\begin{picture}(4824,597)(0,-10)
\put(462,336){\blacken\ellipse{50}{50}}
\put(462,336){\ellipse{50}{50}}
\put(1212,336){\blacken\ellipse{50}{50}}
\put(1212,336){\ellipse{50}{50}}
\put(1962,336){\blacken\ellipse{50}{50}}
\put(1962,336){\ellipse{50}{50}}
\put(3612,336){\blacken\ellipse{50}{50}}
\put(3612,336){\ellipse{50}{50}}
\put(4362,336){\blacken\ellipse{50}{50}}
\put(4362,336){\ellipse{50}{50}}
\put(4587,261){\blacken\ellipse{50}{50}}
\put(4587,261){\ellipse{50}{50}}
\path(12,186)(912,486)
\path(762,486)(1662,186)
\path(1512,186)(2412,486)
\path(3162,186)(4062,486)
\path(3162,186)(4062,486)
\path(3912,486)(4812,186)
\path(3912,486)(4812,186)
\put(462,36){\makebox(0,0)[lb]{\smash{{{\SetFigFont{8}{9.6}{\rmdefault}{\mddefault}{\updefault}$T_1$}}}}}
\put(1062,36){\makebox(0,0)[lb]{\smash{{{\SetFigFont{8}{9.6}{\rmdefault}{\mddefault}{\updefault}$T_2$}}}}}
\put(1962,36){\makebox(0,0)[lb]{\smash{{{\SetFigFont{8}{9.6}{\rmdefault}{\mddefault}{\updefault}$T_3$}}}}}
\put(387,486){\makebox(0,0)[lb]{\smash{{{\SetFigFont{8}{9.6}{\rmdefault}{\mddefault}{\updefault}$p_1$}}}}}
\put(1212,486){\makebox(0,0)[lb]{\smash{{{\SetFigFont{8}{9.6}{\rmdefault}{\mddefault}{\updefault}$p_2$}}}}}
\put(1737,486){\makebox(0,0)[lb]{\smash{{{\SetFigFont{8}{9.6}{\rmdefault}{\mddefault}{\updefault}$p_3$}}}}}
\put(2637,336){\makebox(0,0)[lb]{\smash{{{\SetFigFont{8}{9.6}{\rmdefault}{\mddefault}{\updefault}$\cdots$}}}}}
\put(4287,486){\makebox(0,0)[lb]{\smash{{{\SetFigFont{8}{9.6}{\rmdefault}{\mddefault}{\updefault}$p_{r-j}$}}}}}
\put(4737,261){\makebox(0,0)[lb]{\smash{{{\SetFigFont{8}{9.6}{\rmdefault}{\mddefault}{\updefault}$\infty$}}}}}
\put(3387,486){\makebox(0,0)[lb]{\smash{{{\SetFigFont{8}{9.6}{\rmdefault}{\mddefault}{\updefault}$p_{r-j-1}$}}}}}
\put(3312,36){\makebox(0,0)[lb]{\smash{{{\SetFigFont{8}{9.6}{\rmdefault}{\mddefault}{\updefault}$T_{r-j-1}$}}}}}
\put(4212,36){\makebox(0,0)[lb]{\smash{{{\SetFigFont{8}{9.6}{\rmdefault}{\mddefault}{\updefault}$T_{r-j}$}}}}}
\end{picture}
}
\end{center}
\caption{Breaking the Hurwitz class by breaking the target \lremind{oscar}}
\label{oscar}
\end{figure}
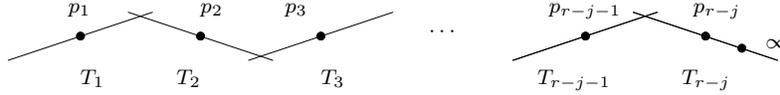

Hence it suffices to show that any one
of the terms of form \eqref{grover} lies in $I_{\geq 2g-2+n-j}$.
Observe that the number of genus $0$ components on the stabilization of
a prestable curve is at least the number of genus $0$ components with at
least 3 special points, minus the number of genus $0$ components with 1
special point.  Then by Lemma~\ref{ilem}, the contribution of any one
of these terms lies in $I_{\geq Q}$, where $Q$ is $r-j-1$ (one for each
$\HH^{g^s, \al^s, \be^s}$, $2 \leq s \leq r-j$) minus at most $d-1$
(for each genus $0$ curve in the preimage of $T_1$ with one special
point; there are at most $d-1$ such components of the preimage of $T_1$).
As $r=d+n+2g-2$ (by the Riemann-Hurwitz formula), the result follows.
\epf

\bpoint{Proof of \thmstar}
The proof of \thmstar{} is now a straightforward generalization
of that of \cite{socle}.

\epoint{Reduction to the case of monomials in $\psi$-classes}
The statement of \thmstar{} behaves well with respect to the
natural morphisms: 
Pushing forward by the forgetful morphism~\eqref{forgetful}
decreases the codimension by $1$, and decreases the 
number of genus $0$ components by at most $1$.  
Pushing forward the
class $a \boxtimes b$ by the gluing morphism \eqref{gluing1}
involves adding both the genera and the numbers of genus $0$ components;
the codimension is the sum of the old codimensions plus $1$.
Pushing forward by the gluing morphism \eqref{gluing2} preserves
number of genus $0$ components, and increases both the genus and codimension by
$1$.  Hence it suffices
to prove \thmstar {} for monomials in the $\psi$-classes,
i.e.\ that
any dimension $j$ monomial in $\psi$-classes lies in $I_{ \geq 2g-2+n-j}$.

\epoint{Localization} We proceed by induction on $\dim \cmbar_{g,n}$.
Apply relative virtual localization (Corollary~\ref{bert}) to the
component of $\cmbar_{g,\al}\left( \proj^1 \right)$ where the source
curve is required to be {\em connected}.  (We could just as well work
on the whole space, but this restriction will simplify our exposition
by avoiding the nonexistent moduli spaces $\M_{0,2}$ and $\M_{0,1}$.)
We cap with $L_{p_1}$, \dots, $L_{p_{r-j}}$ to calculate
$\HH^{g,\al}_{j, \conn}$, defined to be the restriction of
$\HH^{g,\al}_j$ to the component $\cmbar_{g,n}$ of $\cmt$.  Choose
weights on $L_{p_1}$, \dots, $L_{p_{r-j}}$ corresponding to requiring
the fixed branch points points $p_1$, \dots, $p_{r-j}$ to map to $0$.
One fixed locus corresponds to all remaining branch points also
mapping to $0$: the simple fixed locus.  This locus corresponds to an
$n$-pointed genus $g$ curve mapping to $0$, glued to $n$ trivial
covers of $\proj^1$.  By relative virtual localization, the
contribution of this locus is
$$
r! 
\left(  \prod_{i=1}^n \frac {\al_i^{\al_i}}{\al_i!} \right)
\left[  \frac { 1 - \la_1 + \cdots \pm \la_g } { \prod_i ( 1 - \al_i \psi_i)} 
\right]_j
= r! \left(  \prod_{i=1}^n \frac {\al_i^{\al_i}}{\al_i!}  \right) P^g_j(\al_1, \dots, \al_n).
$$
($[x ]_j$ is the dimension $j$ component of $x$)
where $P^g_j$ is a (Chow-valued)
polynomial in $n$ variables whose coefficients include all dimension
$j$ monomials in $\psi$-classes.

We conclude by showing that the contributions from the remaining fixed
loci lie in $I_{\geq 2g-2+n-j}$.  Then as $\HH^{g,\al}_j \in I_{\geq
  2g-2+n-j}$ as well, $P^g_j(\al_1, \dots, \al_n) \in I_{\geq
  2g-2+n-j}$.  As the coefficients of a polynomial lie in the
$\Q$-span of the value of the polynomial evaluated at sufficiently
many lattice points, \thmstar{} follows.

Consider another (necessarily composite) fixed locus, such as the 
one depicted in  Figure~\ref{crazy}.  Suppose
that the preimage of $0$ contains (i) a curve (not necessarily
connected) of genus $g_0$ with $n_0$ marked points, (ii) $N$ nodes,
and (iii) $a''$ isolated smooth points.  Over the non-rigid part of the
target (``over $\infty$''), suppose that there are $a$ trivial covers,
and the remaining part of the cover corresponds
to a connected component of $\cmbar_{g_{\infty}, \be, \al'
}\left( \proj^1 \right)_{\unrigid}$, where $\al'$ is the
partition over $\infty$, and $\be$ is the partition corresponding to
the node of the target.

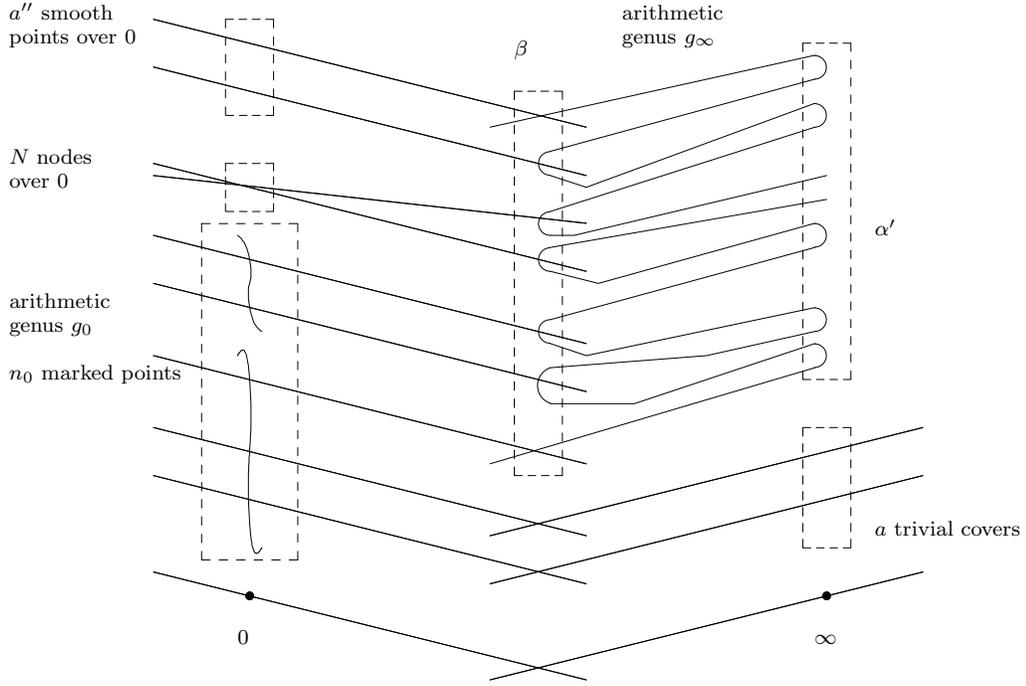
\begin{figure}
\begin{center}
\setlength{\unitlength}{0.00083333in}
\begingroup\makeatletter\ifx\SetFigFont\undefined%
\gdef\SetFigFont#1#2#3#4#5{%
  \reset@font\fontsize{#1}{#2pt}%
  \fontfamily{#3}\fontseries{#4}\fontshape{#5}%
  \selectfont}%
\fi\endgroup%
{\renewcommand{\dashlinestretch}{30}
\begin{picture}(5712,4242)(0,-10)
\put(5025.000,2037.000){\arc{150.000}{4.7124}{7.8540}}
\put(5025.000,3837.000){\arc{150.000}{4.7124}{7.8540}}
\put(5025.000,3537.000){\arc{150.000}{4.7124}{7.8540}}
\put(5025.000,2787.000){\arc{150.000}{4.7124}{7.8540}}
\put(3375.000,3237.000){\arc{150.000}{1.5708}{4.7124}}
\put(3375.000,2862.000){\arc{150.000}{1.5708}{4.7124}}
\put(3375.000,2637.000){\arc{150.000}{1.5708}{4.7124}}
\put(3375.000,2187.000){\arc{150.000}{1.5708}{4.7124}}
\put(3412.500,1849.500){\arc{237.171}{1.8925}{4.3906}}
\put(5025.000,2262.000){\arc{150.000}{4.7124}{7.8540}}
\put(5100,537){\blacken\ellipse{50}{50}}
\put(5100,537){\ellipse{50}{50}}
\put(1500,537){\blacken\ellipse{50}{50}}
\put(1500,537){\ellipse{50}{50}}
\path(900,687)(3600,12)
\path(900,687)(3600,12)
\path(3000,12)(5700,687)
\path(3000,12)(5700,687)
\path(900,1287)(3600,612)
\path(900,1287)(3600,612)
\path(3000,612)(5700,1287)
\path(3000,612)(5700,1287)
\path(3000,912)(5700,1587)
\path(3000,912)(5700,1587)
\path(900,1587)(3600,912)
\path(900,1587)(3600,912)
\path(900,2037)(3600,1362)
\path(900,2037)(3600,1362)
\path(900,2487)(3600,1812)
\path(900,2487)(3600,1812)
\path(900,2787)(3600,2112)
\path(900,2787)(3600,2112)
\path(900,3237)(3600,2562)
\path(900,3237)(3600,2562)
\path(900,3162)(3600,2862)
\path(900,3162)(3600,2862)
\path(900,4137)(3600,3462)
\path(900,4137)(3600,3462)
\path(900,3837)(3600,3162)
\path(900,3837)(3600,3162)
\path(3000,1362)(5025,1962)
\path(3000,3462)(5025,3912)
\path(3375,1737)(3900,1737)(5025,2112)
\path(3375,1962)(4350,2037)(5025,2187)
\path(3375,2112)(3600,2037)(5025,2337)
\path(3375,2262)(5025,2712)
\path(5025,2862)(3675,2487)(3375,2562)
\path(3375,2712)(5100,3012)
\path(3375,2937)(5025,3462)
\path(5025,3612)(3600,3087)(3375,3162)
\path(3375,3312)(5025,3762)
\path(3375,2787)(3525,2787)(5100,3162)
\dashline{60.000}(4950,3987)(5250,3987)(5250,1887)
	(4950,1887)(4950,3987)
\dashline{60.000}(4950,1587)(5250,1587)(5250,837)
	(4950,837)(4950,1587)
\dashline{60.000}(3150,3687)(3450,3687)(3450,1287)
	(3150,1287)(3150,3687)
\dashline{60.000}(1200,2862)(1800,2862)(1800,762)
	(1200,762)(1200,2862)
\dashline{60.000}(1350,4137)(1650,4137)(1650,3537)
	(1350,3537)(1350,4137)
\dashline{60.000}(1350,3237)(1650,3237)(1650,2937)
	(1350,2937)(1350,3237)
\path(1425,2787)(1426,2786)(1428,2785)
	(1432,2782)(1437,2777)(1444,2771)
	(1451,2764)(1458,2754)(1466,2743)
	(1473,2730)(1480,2713)(1487,2693)
	(1494,2667)(1500,2637)(1504,2608)
	(1507,2581)(1508,2557)(1508,2537)
	(1507,2520)(1505,2507)(1503,2497)
	(1500,2487)(1497,2477)(1495,2467)
	(1493,2454)(1492,2437)(1492,2417)
	(1493,2393)(1496,2366)(1500,2337)
	(1506,2307)(1513,2281)(1520,2261)
	(1527,2244)(1534,2231)(1542,2220)
	(1549,2210)(1556,2203)(1563,2197)
	(1568,2192)(1572,2189)(1574,2188)(1575,2187)
\path(1425,2037)(1425,2038)(1426,2038)
	(1426,2039)(1427,2040)(1427,2041)
	(1428,2043)(1429,2045)(1431,2047)
	(1432,2049)(1434,2052)(1436,2054)
	(1437,2057)(1439,2059)(1442,2062)
	(1444,2064)(1446,2066)(1448,2068)
	(1451,2070)(1453,2071)(1456,2072)
	(1458,2072)(1461,2071)(1463,2070)
	(1466,2068)(1468,2066)(1470,2062)
	(1473,2057)(1475,2052)(1478,2045)
	(1480,2037)(1483,2027)(1485,2016)
	(1487,2003)(1490,1988)(1492,1972)
	(1494,1954)(1496,1933)(1498,1911)
	(1500,1887)(1502,1851)(1504,1812)
	(1506,1774)(1507,1736)(1508,1699)
	(1508,1665)(1508,1632)(1508,1602)
	(1507,1574)(1506,1548)(1505,1523)
	(1504,1500)(1503,1479)(1501,1458)
	(1500,1437)(1499,1416)(1497,1395)
	(1496,1374)(1495,1351)(1494,1326)
	(1493,1300)(1492,1272)(1492,1242)
	(1492,1209)(1492,1175)(1493,1138)
	(1494,1100)(1496,1062)(1498,1023)
	(1500,987)(1502,963)(1504,941)
	(1506,920)(1508,902)(1510,886)
	(1513,871)(1515,858)(1517,847)
	(1520,837)(1522,829)(1525,822)
	(1527,817)(1530,812)(1532,808)
	(1534,806)(1537,804)(1539,803)
	(1542,802)(1544,802)(1547,803)
	(1549,804)(1552,806)(1554,808)
	(1556,810)(1558,812)(1561,815)
	(1563,817)(1564,820)(1566,822)
	(1568,825)(1569,827)(1571,829)
	(1572,831)(1573,833)(1573,834)
	(1574,835)(1574,836)(1575,836)(1575,837)
\put(3150,3912){\makebox(0,0)[lb]{\smash{{{\SetFigFont{8}{9.6}{\rmdefault}{\mddefault}{\updefault}$\beta$}}}}}
\put(5025,237){\makebox(0,0)[lb]{\smash{{{\SetFigFont{8}{9.6}{\rmdefault}{\mddefault}{\updefault}$\infty$}}}}}
\put(1425,237){\makebox(0,0)[lb]{\smash{{{\SetFigFont{8}{9.6}{\rmdefault}{\mddefault}{\updefault}$0$}}}}}
\put(5400,2787){\makebox(0,0)[lb]{\smash{{{\SetFigFont{8}{9.6}{\rmdefault}{\mddefault}{\updefault}$\alpha'$}}}}}
\put(0,2187){\makebox(0,0)[lb]{\smash{{{\SetFigFont{8}{9.6}{\rmdefault}{\mddefault}{\updefault}genus $g_0$}}}}}
\put(0,1887){\makebox(0,0)[lb]{\smash{{{\SetFigFont{8}{9.6}{\rmdefault}{\mddefault}{\updefault}$n_0$ marked points}}}}}
\put(0,2337){\makebox(0,0)[lb]{\smash{{{\SetFigFont{8}{9.6}{\rmdefault}{\mddefault}{\updefault}arithmetic}}}}}
\put(0,3237){\makebox(0,0)[lb]{\smash{{{\SetFigFont{8}{9.6}{\rmdefault}{\mddefault}{\updefault}$N$ nodes}}}}}
\put(0,3087){\makebox(0,0)[lb]{\smash{{{\SetFigFont{8}{9.6}{\rmdefault}{\mddefault}{\updefault}over $0$}}}}}
\put(5400,912){\makebox(0,0)[lb]{\smash{{{\SetFigFont{8}{9.6}{\rmdefault}{\mddefault}{\updefault}$a$ trivial covers}}}}}
\put(3825,3987){\makebox(0,0)[lb]{\smash{{{\SetFigFont{8}{9.6}{\rmdefault}{\mddefault}{\updefault}genus $g_{\infty}$}}}}}
\put(3825,4137){\makebox(0,0)[lb]{\smash{{{\SetFigFont{8}{9.6}{\rmdefault}{\mddefault}{\updefault}arithmetic}}}}}
\put(0,3987){\makebox(0,0)[lb]{\smash{{{\SetFigFont{8}{9.6}{\rmdefault}{\mddefault}{\updefault}points over $0$}}}}}
\put(0,4137){\makebox(0,0)[lb]{\smash{{{\SetFigFont{8}{9.6}{\rmdefault}{\mddefault}{\updefault}$a''$ smooth}}}}}
\end{picture}
}
\end{center}
\caption{An example of a composite fixed locus \lremind{crazy}}
\label{crazy}
\end{figure}

By relative virtual
localization, the contribution of this fixed locus is a 
linear combination of terms of the form
$gl({\mathbf c}_{0} \boxtimes {\mathbf c}_{\infty})$ where 
$$
{\mathbf c}_{\infty} = 
\pi_* \left( \psi^k \cap \left[ \cmbar_{g_{\infty}, \be, \al'}\left( \proj^1 \right)_{\unrigid} \right]^{\vir} \right)  \in 
A_{\vdim -k}\left(
\cmbar_{g_{\infty}, \left| \be \right| + \left| \al' \right|}
\right) 
$$
and
$${\mathbf c}_0 \in R_{j-(\vdim-k)}\left( \cmt_{g_0, n_0} \right)$$
where $\vdim = \vdim
(
\cmbar_{g_{\infty}, \be, \al'}\left( \proj^1 \right)_{\unrigid} 
)  = 
 \left| \be \right| +
\left| \al' \right| + 2 g_{\infty}-3$.

By Lemma~\ref{l3}, ${\mathbf c}_{\infty} \in I_{\geq k+1}$, and by the
inductive hypothesis and Remark~\ref{ernie}, ${\mathbf{c}}_0 \in I_{\geq 2 g_0 - 2 + n_0 - ( j-(\vdim-k))}$.
Thus we wish to show that
\begin{equation} \label{toprove}
(k+1) + (2g_0-2+n_0-j+ \left| \be \right| + \left| \al' \right| + 2 g_{\infty} - 3 - k) \geq 2g-2+n-j.
\end{equation}
By comparing the genus of the source curve to that of its components, we have
$$
(g_0-1) + (g_{\infty} - 1 ) - N  - a'' +  \left| \be \right| = g-1.
$$
By counting preimages of the node of the target in two ways,
$$
n_0 + 2N + a'' = \left| \be \right| + a.
$$
Finally, the number of preimages of $\infty$ is
$$
a + \left| \al' \right| = n.
$$
By adding twice the first equation to the other two, we have
$$
\left( (k+1) + \left( 2g_0-2+n_0-j+ \left| \be \right| 
  + \left| \al' \right| + 2 g_{\infty} - 3 - k \right)
\right)  -a'' = 2g-2+n-j  
$$
from which \eqref{toprove} follows.
Thus 
$gl({\mathbf c}_0 \boxtimes
{\mathbf c}_{\infty}) \in I_{\geq 2g + \left| \al \right| - 2 - j}$ as desired. 
\epf

\section{Applications of \thmstar}
\label{consequences}\lremind{consequences}We 
give applications of \thmstar{} to prove and extend
various theorems and conjectures.  
We will use the bijection between strata
of $\mgn$ and stable graphs, and the ``cross-ratio'' or WDVV relation
among strata. 
We also will make repeated use of the
following.

\tpoint{Corollary}   {\em \lremind{impcor}\label{impcor}Any $c \in R_j \left( \cmbar_{g,n} \right)$ is the pushforward
(under inclusion) of classes supported on boundary strata  corresponding
to curves with components of geometric genus $g_k$ and $n_k$ special points
(marked points and node-branches),
where
$$
\sum_k (g_k-1+\de_{g_k,0}) +
 \sum_k (2 g_k - 2 + n_k - \de_{g_k,0})
\geq j \geq \sum_k (2 g_k - 2 + n_k - \de_{g_k,0}).
$$} 
(As usual, $\de_{g_k,0}=1$ if $g_k=0$, and $0$ otherwise.)
Note that $2 g_k - 2
+ n_k - \de_{g_k,0} \geq 0$ with equality if and only if $g_k = 0$ and
$n_k=3$.  Also, if $g_k=0$ or $1$ for all $k$, then the contribution
of this stratum to $c$ is a multiple of the fundamental class of the stratum.

\bpf
The left inequality is the dimension of the stratum.
By \thmstar, 
$$\sum \de_{g_k,0} \geq (3g-3+n-j)-g+1 = 2g-2+n-j = \sum(2g_k-2+n_k) - j,$$ 
from which the right inequality follows. \epf

\bpoint{Getzler's conjecture (Ionel's theorem)} In
\cite[Footnote~1]{gconj}, Getzler conjectured that
all degree $g$ monomials in $\psi_1$, \dots, $\psi_n$ vanish on
$\cm_{g,n}$ if $g>0$.  Getzler's conjecture was known in genus $1$
(classically), and genus $2$ (by work of Mumford and Getzler,
see equs. (4) and (5) of \cite{gconj}).  Ionel proved Getzler's
conjecture in cohomology \cite{igetz}.  Her
argument, rewritten in the language of algebraic geometry, should also
prove the conjecture in Chow.

\thmstar {} immediately implies Getzler's conjecture, as well as more: {\em
  all} tautological classes of codimension at least $g$ vanish on
$\cm_{g,n}$ if $g>0$, and in fact they vanish on the larger open set
corresponding corresponding to curves with no rational component.

\bpoint{Poincar\'{e} duality speculations}
The next few applications concern three parallel conjectures or speculations.
Let $\mgnrt$ be the moduli space of ``curves with rational tails'',
curves with a smooth component of genus $g$ (i.e.\
with dual graph with a vertex of genus $g$). 
Let $\mgnct$ be the moduli space of ``curves of compact type'', 
curves with compact Jacobian (i.e.\ with dual graph with no loops).
Hence
$$\mgnrt \subset \mgnct \subset \mgn.$$

\tpoint{Conjecture (Faber, Looijenga, Pandharipande, et al.\ \cite[Sect.~2]{icm})}
{\em The space\label{csq}\lremind{csq} 
$\mgn$ (resp.\ $\mgnct$, $\mgnrt$) ``behaves like''
a complex variety of dimension $D = 3g-3+n$ (resp.\ $2g-3+n$, $g-2+n$).
More precisely,
\begin{itemize}
\item {\rm Socle statement:} $R^i = 0$ for $i>D$, $R^D \cong \Q$, and 
\item {\rm Perfect pairing statement:} for $0 \leq i \leq D$, the natural map $R^i \times R^{D-i} \rightarrow
R^D$ is a perfect pairing.
\end{itemize}}

Although it is not clear if one should expect this strong statement to
be true, this speculation has motivated much interesting work.

The case $\cm_{g,0}=\cm^{rt}_{g,0}$ is part of Faber's conjecture \cite{fconj}.  The
case $\mgn$ was asked by Hain and Looijenga 
\cite[Question~5.5, p.~108]{hl}.  The cases when
$n=0$ were stated in \cite{fpconj}.  The general statement is likely
due to Faber and Pandharipande.

Some evidence for the perfect pairing portion will be given in a
later paper \cite{gvlater}.  We now present proofs of the socle
statements in all three cases. 

\epoint{Socle proof for $\mgn$} \label{soclestable}\lremind{soclestable}This 
argument is essentially the one given in
\cite{socle}, which can be interpreted as an early version of 
\thmstar.  If $i>3g-3+n = \dim \mgn$, $R^i\left( \mgn \right)=0$ for dimensional
reasons.  If $i=3g-3+n$ then by Corollary~\ref{impcor}, $R^i\left( \mgn \right)$ is
generated by the 0-dimensional strata (whose graphs have $2g-2+n$
genus $0$ trivalent vertices).  These strata are rationally equivalent
(by judicious use of the cross-ratio relation, or by observing that
they lie on the image of the rational space $\cmbar_{0,2g+n}$ under
$g$ gluing morphisms \eqref{gluing2}).  They are non-zero
because they have nonzero degree.

\epoint{Socle proof for $\mgnct$} 
\label{soclemgnct}\lremind{soclemgnct}This is the first genus-free evidence for the ``compact type'' conjecture.
As a bonus, this approach produces a natural generator of the socle.

The argument parallels the previous one.  The case $(g,n) = (2,0)$ is
immediate.  For $(g,n) \neq (2,0)$, any
$n$-pointed genus $g$ stable graph with no loops has no more than
$g-2+n$ genus $0$ vertices, and if equality holds, then all genus 0
vertices are trivalent, and the other vertices are genus $1$
``leaves'' (see Figure~\ref{mgnct} for an example),
so in particular Corollary~\ref{impcor} applies.  Any such strata is codimension
$2g-3+n$ (dimension $g$).  Any two such strata are rationally
equivalent by the cross-ratio relation, or by using the rationality of
the space $\cmbar_{0,g+n}$.

\begin{figure}
\begin{center}
\setlength{\unitlength}{0.00083333in}
\begingroup\makeatletter\ifx\SetFigFont\undefined%
\gdef\SetFigFont#1#2#3#4#5{%
  \reset@font\fontsize{#1}{#2pt}%
  \fontfamily{#3}\fontseries{#4}\fontshape{#5}%
  \selectfont}%
\fi\endgroup%
{\renewcommand{\dashlinestretch}{30}
\begin{picture}(3525,1197)(0,-10)
\put(1275,675){\blacken\ellipse{50}{50}}
\put(1275,675){\ellipse{50}{50}}
\put(1575,225){\blacken\ellipse{50}{50}}
\put(1575,225){\ellipse{50}{50}}
\put(2025,525){\blacken\ellipse{50}{50}}
\put(2025,525){\ellipse{50}{50}}
\put(2325,975){\blacken\ellipse{50}{50}}
\put(2325,975){\ellipse{50}{50}}
\put(1275,975){\ellipse{100}{100}}
\put(2775,375){\ellipse{100}{100}}
\put(2625,825){\ellipse{100}{100}}
\put(2025,1125){\ellipse{100}{100}}
\put(2175,75){\ellipse{100}{100}}
\put(2475,225){\blacken\ellipse{50}{50}}
\put(2475,225){\ellipse{50}{50}}
\path(3465,998)(2418,978)
\blacken\path(2537.405,1010.286)(2418.000,978.000)(2538.551,950.297)(2501.985,979.604)(2537.405,1010.286)
\path(3489,323)(2872,374)
\blacken\path(2994.063,394.013)(2872.000,374.000)(2989.121,334.217)(2955.715,367.080)(2994.063,394.013)
\path(2326,967)(2067,1097)
\path(2340,964)(2579,841)
\path(2480,214)(2729,343)
\path(2473,217)(2217,108)
\path(1579,217)(1320,98)
\path(1282,664)(1020,548)
\path(1282,667)(1275,920)
\path(1275,675)(1575,225)(2025,525)(2475,225)
\path(2325,975)(2025,525)
\path(525,790)(859,541)
\blacken\path(744.862,588.671)(859.000,541.000)(780.724,636.775)(791.655,591.206)(744.862,588.671)
\put(900,450){\makebox(0,0)[lb]{\smash{{{\SetFigFont{8}{9.6}{\rmdefault}{\mddefault}{\updefault}1}}}}}
\put(1200,0){\makebox(0,0)[lb]{\smash{{{\SetFigFont{8}{9.6}{\rmdefault}{\mddefault}{\updefault}2}}}}}
\put(3525,975){\makebox(0,0)[lb]{\smash{{{\SetFigFont{8}{9.6}{\rmdefault}{\mddefault}{\updefault}genus $0$}}}}}
\put(3525,300){\makebox(0,0)[lb]{\smash{{{\SetFigFont{8}{9.6}{\rmdefault}{\mddefault}{\updefault}genus $1$}}}}}
\put(0,825){\makebox(0,0)[lb]{\smash{{{\SetFigFont{8}{9.6}{\rmdefault}{\mddefault}{\updefault}marked point}}}}}
\end{picture}
}
\end{center}
\caption{An example of a generator of $R^{2g-3+n}(\cm^{ct}_{5,2})$}
\label{mgnct}
\end{figure}

Hence by \thmstar, if $i>D=2g-3+n$, then any element of $R^i(\mgnct)$
vanishes on $\mgnct$, and any element of $R^D(\mgnct)$ is a multiple
of (any) one of these strata.  This strata is nonzero on $\mgnct$ as
the integral of $\la_g$ over its closure is clearly
$(\int_{\cmbar_{1,1}} \la_1)^g = 1/24^g \neq 0$, and $\lambda_g$
vanishes on $\cmbar_{g,n} - \mgnct$ \cite[Sect.\ 0.4]{fpconj}.

We remark that these curves are related to {flag curves} 
\cite[p.~246]{hm}.

\epoint{Socle proof for $\mgnrt$} \label{soclemgnrt}\lremind{soclemgnrt}Conjecture~\ref{csq} is known for $g=0$ (\cite{keel}, as
$\cmbar_{g,0} = \cm_{g,0}^{rt}$), so we assume $g>0$.  

\tpoint{Proposition}  
{\em  
\begin{enumerate}
\item[(a)]  For $i > g-2+n$, $R^i(\mgnrt) = 0$.
\item[(b)]
For $n>0$, the forgetful morphism $\al:  R^{g-2+n}(\mgnrt) \rightarrow R^{g-1}(\cm_{g,1}^{rt} = \cm_{g,1})$ is an isomorphism.
\end{enumerate}
}

Part (a) for $n=0$ is Looijenga's theorem \cite{lthm}. 
The Proposition should be deducible from the main theorem of \cite{lthm},
although we have not checked the details.

\bpf
(a)  If $n>0$, then $\mgnrt \subset 
\cS_{\leq n-1}$, so 
by \thmstar{} any class in $R^i(\mgnrt)$ is the pushforward
of a class supported on $\cS_{\geq n} \cap \mgnrt = \emptyset$.
By Definition~\ref{tautdef} of the tautological
ring, $R^{a+1}(\cm_{g,1}) \rightarrow R^a(\cm_g= \cm_g^{rt})$ is surjective,
so the case $n=0$ follows.

(b)  Any $[C] \in \cmbar_{0,n+1}$ induces a closed immersion
$\iota_C:  \cm_{g,1} \rightarrow \mgnrt$ via the gluing morphism
\eqref{gluing1}.  As $\cmbar_{0,n+1}$ is rational, the induced map $A_*(\cm_{g,1}) \rightarrow A_*(\mgnrt)$
is independent of $C$; denote it by $\be$.
Let $C_1$, \dots, $C_s$ be the curves corresponding to the 0-dimensional
strata of $\cmbar_{0,n+1}$.  By \thmstar, any $a \in R^{g-2+n}(\mgnrt)$
is the pushforward of a class supported on $\cS_{n-1} \cap \mgnrt = \cup_j \iota_{C_j}(\cm_{g,1})$,
and hence the pushforward of a class $b$ supported on $\iota_{C_1}(\cm_{g,1}) \cong \cm_{g,1}$.  Necessarily $b = \al(a)$ (and thus is tautological).
Hence $\be \circ \al$ is the identity on $R^* \left( \cmbar_{g,n} \right)$.  As
$\al \circ \be$ is clearly the identity on $R^* \left( \cm_{g,1} \right)$, we are done.
\epf

In light of Faber's calculation $\deg \kappa_{g-2} \lambda_g
\lambda_{g-1} \neq 0$ and that $\lambda_g \lambda_{g-1}$ vanishes on
$\cmbar_g - \cm_g$ \cite[Theorem~2]{fconj}, it now follows that
$R^{g-2+n}(\mgnrt) \neq 0$ for all $g,n$.  Finally, $R^{g-2+n}(\mgnrt)
\cong \Q$ follows from $R^{g-1}(\cm_{g,1}) = \langle \psi_1^{g-1}
\rangle$, which is known by Looijenga.  We do not see how to obtain
this result easily using our approach.

\bpoint{Generalizations of Diaz' theorem}\label{diazgen} Following
\cite[p.~412]{lthm}, Diaz' theorem \cite[p.~79]{diaz} 
(there is no
  complete subvariety of $\cm_g$ of dimension greater than
  $g-2$)
can be extended as follows.

\tpoint{Proposition}  {\em 
\begin{enumerate}
\item[(a)] There is no
  complete subvariety of $S_{\leq s}$ of dimension greater than
  $s+g-1$. In particular, there is no complete subvariety of $S_0$ of dimension
greater than $g-1$.
\item[(b)] There is no complete subvariety of $M_{g,n}^{rt}$ of
dimension greater than $g-2+n$.
\item[(c)] There is no complete subvariety of
  $M_{g,n}^{ct}$ of dimension greater than $2g-3+n$. 
\end{enumerate}}

Here roman letters denote the coarse moduli space.  Part (b) is a
trivial generalization of Diaz' theorem, and part (c) follows from
another result of Diaz \cite[Cor.~ p.~80]{diaz} (\cite{fl} notes that
$2g-1$ should be $2g-2$ here); we include the argument only as an
illustration of how these results are also immediately consequences of
\thmstar.  Note that the bound in (a) for $s=0$ is necessarily worse
than the bound of Diaz' theorem, as there is a complete curve in $S_0(
\cmbar_{2,0})$ (e.g.\ $y^2 = x(x-1)(x-2)(x-3)(x-4)(x-t)$).

\bpf Suppose $V$ is a complete substack of $\cS_{\leq s}$ of dimension
$d \geq g-1$.  Let $\al$ be any ample tautological divisor on the
coarse moduli space.  (The coarse moduli space is projective by
Knudsen, and all divisors are tautological by \cite{acharer}.
Alternatively, Cornalba gives a tautological ample class explicitly
\cite[p.\ 13]{cornalba}.)  If $d \geq e \geq g-1$, the class
$\al^e$ is tautological, and by \thmstar {} is rationally equivalent to
a class supported on $\cS_{\geq e-g+1}$.  Hence $\cS_{\geq e-g+1} \cap
V$ contains a nonempty subset of dimension at least $d-e$.  The case
$d=e=s+g$ yields (a). Parts (b) and Part (c) follow from $M_{g,n}^{rt}
\subset S_{\leq n-1}$ ($n \geq 1$) and $M_{g,n}^{ct} \subset S_{\leq
  g-2+n}$, or directly from Sections~\ref{soclemgnrt}
and~\ref{soclemgnct} respectively.  \epf

\bpoint{Universal description of the tautological groups in low
  dimension} \lremind{lowdim}\label{lowdim}The tautological groups are
completely understood in codimension $1$ and $2$ (\cite{acharer} and
\cite{polito} respectively), but even codimension 3 is combinatorially
complicated.  \thmstar {} implies that the tautological groups are
actually more straightforward in low dimension, by providing a
parsimonious systems of generators.  This is by exploiting the
right-hand inequality of Corollary~\ref{impcor}: for small values of
$j$, there are very few possible values of $g_k$ and $n_k$ that can
occur in a boundary stratum.  For notational convenience, we will
imprecisely denote such a stratum (with components of genus $g_k$,
with $n_k$ special points) by $\prod_k \M_{g_k,n_k}$.  We will list the
consequences of this inequality up through dimension $6$.

\noindent
{\em Dimension $0$.} As observed in Section~\ref{soclestable}, 
$R_0\left( \mgn \right)$
is generated by boundary strata of the form $\prod \M_{0,3}$; 
they are all rationally equivalent.  

\noindent
{\em Dimension $1$.}
$R_1\left( \mgn \right)$ is again generated by 
boundary strata (in
which there is one $\cmbar_{0,4}$ or $\cmbar_{1,1}$ factor).  To see
this, note that by the right hand inequality of Corollary~\ref{impcor}
with $j=1$, the only pairs $(g_k,n_k)$ that can occur are 
$(0,3),(0,4),$ or $(1,1)$, and for all but one value of 
$k$, we must have $(g_k,n_k)=(0,3)$.  Since these boundary strata
have dimension $1$, the only one-dimensional classes they can support
are their fundamental classes.

\noindent
{\em Dimension $2$.}
In dimension $2$, Corollary~\ref{impcor} allows for the
possibility of a non-boundary class: a divisor on $\M_{2,0}$.  However,
by \cite[Part~III]{m}, $A^1 \left( \cmbar_{2,0} \right)$ 
is generated by boundary
strata, so we can conclude that 
$R_2\left( \mgn \right)$ is generated by boundary
strata for all $g$ and $n$.

\noindent {\em Dimension $3$.} $R_3\left( \mgn \right)$ is 
generated by boundary strata, together with
classes supported on strata corresponding to a product of 
$\cmbar_{2,1}$ with
copies of $\cmbar_{0,3}$; there
is a single nonboundary generator
corresponding to $\psi_1$ on
$\cmbar_{2,1}$, as 
this is the only nonboundary class in $A^1(\M_{2,1})$.

\noindent
{\em Dimension $4$.} $R_4\left( \M_{3,0} \right)$ was computed 
in \cite{fgenus3}.  For other $g$, $n$, 
$R_4 \left( \mgn \right)$ is generated by boundary classes,
as well as classes  of the form $\psi_1$ on either
$$\cmbar_{2,2}\times \prod \cmbar_{0,3},  \quad
\M_{2,1} \times \M_{0,4} \times \prod \M_{0,3}, \quad  \text{or}
\quad  \M_{2,1} \times \M_{1,1} \times \prod \M_{0,3}$$ 
with the $\psi$ class 
on the genus $2$ factor.
The last two cases require the facts
that 
$A^1 \left( \cmbar_{2,1} \times \cmbar_{0,4} \right) 
\cong A^1 \left( \cmbar_{2,1} \right)
\oplus A^1 \left( \cmbar_{0,4} \right)$ and 
$A^1 \left( \cmbar_{2,1} \times \cmbar_{1,1} \right) 
\cong A^1 \left( \cmbar_{2,1} \right)
\oplus A^1 \left( \cmbar_{1,1} \right),$
  which follow from $\cmbar_{0,4} \cong \overline{M}_{1,1} \cong
  \proj^1$.

\noindent 
{\em Dimension $5$.} The story is essentially the same:  the cases
$(g,n) = (3,0)$ and $(3,1)$ follow from \cite{fgenus3, fgenus4}; in other cases
we have boundary strata, and classes arising from $\psi_1$ on 
$\cmbar_{2,j}$ ($1 \leq j \leq 3$).

\noindent
{\em Dimension $6$.} Here, the first uncertainty arises: one
possibility is a codimension $2$ class on $\cmbar_{3,2} \times \prod
\cmbar_{0,3}$.  This class  a priori needn't be
tautological on $\cmbar_{3,2}$, and $A^2 \left( \cmbar_{3,2} \right)$ has not yet been
computed.  (There is no technical obstruction to this computation.
Presumably $R^2 \left( \cmbar_{3,2} \right) = A^2 \left( \cmbar_{3,2}
\right)$.)


One would expect that the relations between these generators of
the low dimensional tautological groups will be generated
by relations in low genus.  For example, they are generated by the
cross-ratio relations in dimension $0$ and $1$ (Section~\ref{soclestable}
and \cite{gvlater} respectively), and they are likely generated by the
cross-ratio relations and Getzler's relation \cite{m14} in dimension
$2$.

\bpoint{Consequences in low genus} Many interesting consequences of
\thmstar {}  in low genus already follow from
Getzler's conjecture and known facts about the moduli
space.  For example, in genus $1$, Corollary~\ref{impcor} immediately
shows that the tautological groups are generated by boundary strata;
but this follows from the classical formula $\psi_1 = \de_0/12$.  In
genus 2, Corollary~\ref{impcor} implies that the tautological groups
are generated by boundary strata, and the divisor corresponding to
$\left. \psi_1 \right|_{\cmbar_{2,n}}$ on the strata corresponding to gluing
$\cmbar_{2,n}$ to copies of $\cmbar_{0,n_i}$ and $\cmbar_{1,n_j}$.
This follows from Mumford's formula for $\psi_1^2$ on $\cmbar_{2,1}$
and Getzler's formula for $\psi_1 \psi_2$ on $\cmbar_{2,2}$
(equs.\ (4) and (5) of \cite{gconj}).

Conversely, \thmstar {} (and even Getzler's
conjecture) provides a straightforward proof of both Mumford's and
Getzler's formulas.  For example, from \thmstar, on $\cmbar_{2,1}$,
$\psi_1^2$ is the sum of boundary strata, and five test families
determine the multiplicities with which they occur.

} 

\end{document}